\documentclass[10pt]{article} 
\usepackage{geometry} 
\usepackage[utf8]{inputenc}
\usepackage{titlesec}
\usepackage{amssymb}
\usepackage{amsmath}
\usepackage{amsfonts}
\usepackage{amsthm}
\usepackage{stmaryrd}
\usepackage{mathrsfs}
\usepackage{multicol}
\usepackage{algorithmic}
\usepackage{tikz} 
 \usetikzlibrary{matrix,decorations.pathreplacing,calc}
\usetikzlibrary{decorations.pathreplacing}
\usetikzlibrary{patterns}
\usepackage{graphics}
\usepackage{listings}
\usepackage{xcolor}
\usepackage{pgfplotstable} 
\usepackage{booktabs} 
\usepackage{subcaption}
\usepackage{bm} 
\usepackage{tikz}
\usetikzlibrary{arrows.meta, positioning, shapes.geometric}
\usetikzlibrary{fit,positioning}
\usetikzlibrary{positioning}

\usepackage{xcolor}

\usepackage{algorithm, algorithmic}

\usepackage{amsopn}
\usepackage{mathtools}
\usepackage{pgfplots}
\usepackage{pgfplotstable}
\usepackage{booktabs}
\usepackage{cite}
\definecolor{myrefcolor}{rgb}{0.0, 0.0, 0.5} 
\usepackage[pagebackref=false,colorlinks,linkcolor=blue,citecolor=myrefcolor]{hyperref}

\usepackage{tikz}
\usetikzlibrary{matrix,decorations.pathreplacing,calc}

\usepackage{caption}
\usepackage{subcaption}
\usepackage{todonotes}

\newtheoremstyle{mystyle}
{3pt} 
{3pt} 
{\itshape} 
{} 
{\bfseries} 
{.} 
{ } 
{} 

\newtheoremstyle{regularstyle}
{3pt} 
{3pt} 
{} 
{} 
{\bfseries} 
{.} 
{ } 
{} 

\theoremstyle{mystyle}
\newtheorem{theorem}{Theorem}

\newtheorem{lemma}[theorem]{Lemma}
\newtheorem{definition}[theorem]{Definition}
\newtheorem{proposition}[theorem]{Proposition}

\theoremstyle{regularstyle}
\newtheorem{remark}[theorem]{Remark}
\newtheorem{example}[theorem]{Example}

\numberwithin{equation}{section}
\numberwithin{theorem}{section} 

\usepackage{fancyhdr}
\title{\bf \Large Exploiting the nonzero diagonal pattern in matrix function computations}

\author{\normalsize{Majed Hamadi$^1$, Nezam Mahdavi-Amiri$^1$ and Marcel Schweitzer$^2$
}\\
\footnotesize{\em $^1$Department of Mathematical Sciences, Sharif University of Technology, Tehran, Iran.}\\\footnotesize{\em Emails: }%
\footnotesize{\em\texttt{majed.hamadi77@sharif.edu}},\\
\footnotesize{\em \texttt{nezamm@sharif.edu}}
\\
\footnotesize{\em $^2$School of Mathematics and Natural Sciences, Bergische Universität Wuppertal, Wuppertal, Germany.}\\\footnotesize{\em Email: \texttt{marcel@uni-wuppertal.de}}
}

\graphicspath{{Images/}} 
\date{}

\definecolor{myblue}{rgb}{0.1, 0.2, 0.5}
\titleformat{\section}[runin]
{\normalfont\large\bfseries}{\thesection}{1em}{}[.~\ ]
\titlespacing{\section}{0pt}{1.5ex plus 0.5ex minus 0.2ex}{*0}

\titleformat{\subsection}[runin]
{\normalfont\normalsize\bfseries}{\thesubsection}{1em}{}[.~\ ]
\titlespacing{\subsection}{0pt}{1.5ex plus 0.5ex minus 0.2ex}{*0}

\titleformat{\subsubsection}[runin]
{\normalfont\normalsize\bfseries}{\thesubsubsection}{1em}{}[.~\ ]
\titlespacing{\subsubsection}{0pt}{1.5ex plus 0.5ex minus 0.2ex}{*0}
\newcommand{\dist}{\mathrm{d}}
\begin{document}
\maketitle
\vspace*{-0.7cm}
\begin{abstract}
We consider the task of approximating a matrix function $f(A)$, where $A$ is a matrix in which only a relatively small number of (not necessarily consecutive) sub- and superdiagonals contain nonzero entries. Approximating $f$ by a low-degree polynomial $p$ allows us to obtain sparse approximations to $f(A)$, which one can efficiently work with (while, in general, $f(A)$ is a dense matrix, even when $A$ is sparse). Our approach is based on carefully inspecting the locations where nonzeros can occur in $p(A)$, and identifying the entries in $A$ that influence them. In particular, we illustrate how this approach can be used for efficiently approximating the trace of $f(A)$ and identify how this approach is related to established (stochastic) probing methods for trace estimation. Another application area in which our approach works particularly well is the computation of functions of Toeplitz matrices. Here, studying the sparsity pattern of $p(A)$ allows us to reduce the computation of the whole matrix polynomial to that of a single small-scale submatrix, yielding an algorithm that scales exceptionally well to large problem sizes.
\end{abstract}

\vspace*{0.1cm}
{\footnotesize

}
\vspace*{0.1cm}

\section{Introduction}
Computing matrix functions such as matrix exponentials has recently attracted much attention, as they appear in many important applications, such as the numerical solution of partial differential equations \cite{Dur,Lee}, electronic structure calculations \cite{Bekas-Kokiopoulou-Saad,Goedecker}, and social network analysis \cite{Estrada-Higham1}; see \cite{NHigham} for a comprehensive overview of the topic.

For $A \in \mathbb{C}^{n \times n}$ and $f$ an analytic function on and inside a closed contour $\Gamma$ that encloses the set of eigenvalues of $A$, the corresponding matrix function $f(A)$ can be represented by a contour integral,
\[
f(A) := \dfrac{1}{2\pi \textbf{i}}\int_{\Gamma} f(z) (zI-A)^{-1} \mathrm{d}z;
\]
see \cite{NHigham}. One major difficulty arising in computations involving large-scale matrix functions is that, in general, $f(A)$ is a dense matrix even when $A$ is very sparse. This often makes it infeasible to explicitly construct $f(A)$, due to both high computational complexity and enormous storage demands.

Let us mention that in some applications, not $f(A)$ itself is required, but rather its action on a vector, $f(A)b$, or the value of a quadratic form, $b^Tf(A)b$. In this case, it is not necessary---and not recommended---to explicitly compute $f(A)$ first. Instead, one aims to directly approximate the sought-after quantity by an iterative method, typically by a polynomial~\cite{drkn89,Saad1992} or rational~\cite{Guettel2010} Krylov subspace method. 

Here, we consider the situation that one is indeed explicitly interested in $f(A)$ itself (or its trace). One situation in which this \emph{is} feasible is when $A$ is sparse (or banded) and the function $f$ can be accurately approximated by a low-degree polynomial $p$ on the numerical range $\mathcal{W}(A)$ of $A$; i.e., $f(z) \approx p(z)$ for $z \in \mathcal{W}(A)$. In this case, one can approximate $f(A) \approx p(A)$. This is beneficial because the sparsity pattern of $p(A)$ grows in a controlled fashion with increasing degree of $p$. If $p$ has degree $k$, then $p(A)$ can only have nonzeros in positions $(i,j)$ where $i$ and $j$ have a distance of at most $k$ in the graph of $A$; see, e.g.,~\cite{Frommer-Schimmel-Schweitzer}. As a particular---and widely studied~\cite{Ben_Boit, Benzi-Golub, Ben_Raz,Demko}---special case, if $A$ is banded with bandwidth $\beta$, then $p(A)$ is banded with bandwidth $k\beta$. Thus, as long as $k$ is moderate, it will be computationally feasible to work with $p(A)$.

Here, we investigate a situation that is slightly more general than the banded case: we consider matrices in which only a small number of sub- or superdiagonals contain nonzero elements, but these do not necessarily appear consecutively in the matrix. Of course, such matrices can be treated as banded matrices (where the band is not densely populated), but it is computationally beneficial to incorporate precise knowledge of the exact locations of the nonzero diagonals. We derive an efficient computational procedure for identifying the locations of nonzero elements in $p(A)$, which we then employ as a basic building block for several different algorithms dealing with matrix functions. In particular, we consider the task of approximating $f(A)$ itself (as outlined above) but also the task of only approximating certain entries of $f(A)$ or its trace. 

It turns out that our approach leads to a particularly attractive algorithm when $A$ is a Toeplitz matrix. In this case, it is possible to reduce the computation of $p(A)$ to the computation of a function of a single much smaller submatrix $A$.

\subsection{Outline}\label{subsec:outline}

The remainder of our work is organized as follows. We end this section by introducing some notations. In Section \ref{sec2}, we present a few basic results as foundations for the algorithms we develop and recall some important needed results from the literature. In Section \ref{sec3}, we explain our methodology for identifying positions of nonzero elements in $p(A)$. Building upon this, we identify which submatrices of $A$ play a role in computing the individual entries $[p(A)]_{ij}$ in Section~\ref{sec4}. This forms the basis for algorithms that can be used for approximating the trace of $f(A)$  which we discuss in Section~\ref{sec5}. We then explain how our methodology can be used for efficiently computing functions of Toeplitz matrices in Section~\ref{sec6}. Numerical experiments for illustrating the behavior of the proposed algorithms are presented in Section \ref{sec7} before drawing our conclusions in Section \ref{sec8}. We give some guidelines for how to find a suitable polynomial approximation degree in Appendix~\ref{appendix}.

\subsection{Notation}\label{subsec:notation}
We denote the entry in position $(i,j)$ of a matrix as $[A]_{ij}$. The $r$th diagonal of a matrix $A$, denoted as $\mathtt{diag}(A,r)$, consists of all the matrix elements $[A]_{ij}$ with $j-i = r$. Thus, $\mathtt{diag}(A,0)$ denotes the main diagonal, while $r < 0$ corresponds to elements below the main diagonal and $r > 0$ corresponds to elements above the diagonal. For later use, we denote the set of indices corresponding to all possible diagonals as $\mathcal{D}_n:=\lbrace -n+1,\ldots,n-1\rbrace$. We further denote the set of indices of nonzero diagonals of $A$ as $\mathrm{ND}(A) := \left\{ r\in \mathcal{D}_n:~ \mathtt{diag}(A,r) \neq \mathbf{0} \right\}$. Note that throughout the paper, we assume that $\mathrm{ND}(A)$ is known.

For $I, J \subseteq \{ 1, 2,\ldots, n\}$, we denote by $A[I, J]$ the submatrix of $A$ corresponding to the row indices $I$ and column indices $J$. Further, we define the order function $\mathbb{O}: I \rightarrow \{1,\ldots, |I|\} $ on the index set $I$, where $| I |$ denotes the number of elements in the set $I$, and for $j\in I$, $\mathbb{O}(j) = |\{ i\in I: i\leq j\}|$.\footnote{For example, if $I = \{3,5,9\}$, then $\mathbb{O}(3) = 1$, $\mathbb{O}(5) = 2$ and $\mathbb{O}(9) = 3$.} For $I$ and $J$ subsets of integer numbers, we define $I + J: = \{ i+j : i \in I,~ j \in J \} $ and $ I - J := \{ i-j : i \in I,~ j \in J \} $. We also define the sum and difference between an integer number $i$ and a set $J$ analogously; i.e., $ i + J := \{ i + j : j \in J \} $ and $i - J := \{ i - J : j \in J \}$.

The set of all polynomials of degree at most $k$ is denoted as $\Pi_k$ and the Frobenius matrix norm is denoted by $\Vert \cdot \Vert_F$.

\section{Basic material -- polynomials of sparse and banded matrices}\label{sec2}
In this section, we present a collection of key results concerning matrices and polynomials, which form the foundational framework for our methods and results. 

An established way to obtain scalable algorithms in applications that require working with a matrix function $f(A)$---and not just its actions on vectors---is to exploit the fact that low-degree polynomials in $A$ are sparse when $A$ is sparse~\cite{Benzi-Boito-Razouk, Ben_Raz}. Thus, replacing $f(A) \approx p(A)$ with $p \in \Pi_k$ yields a sparse approximation of $f(A)$ as long as $k$ is moderate.

This was first observed and exploited for \emph{banded} matrices~\cite{Demko,Ben_Raz}: if $A$ is $\beta$-banded (i.e., $[A]_{ij} = 0$, whenever $|i-j| >  \beta$), then $p_k(A)$ is $k\beta$-banded. To extend this to arbitrary sparse matrices $A \in \mathbb{C}^{n \times n}$, one considers the \emph{graph of $A$}.

\begin{definition}\label{def:GA}
The (directed) graph $G(A) = (V,E)$  of a sparse matrix $A \in \mathbb{C}^{n \times n}$ has nodes $V = \{1,\dots,n\}$ and edges $E = \{(i,j): a_{ij} \neq 0, i\neq j\}.$

We denote the geodesic distance, that is, the number of edges on the shortest path, from node $i$ to node $j$ in $G(A)$ by $\dist(i,j)$, where $\dist(i,j) = 0$ if $i=j$.
\end{definition}

One easily shows that $A^k$ only contains nonzeros in positions $(i,j)$ with $\dist(i,j) \leq k$; i.e., nonzeros spread along paths in $G(A)$. If $G(A)$ has a large diameter\footnote{The diameter of $G(A)$ is the largest geodesic distance between any pair of vertices in $G(A)$.}---as is typical for sparse matrices in many applications---$p(A)$ with $p \in \Pi_k$ is also rather sparse for moderate values of $k$.\footnote{If, in contrast, $G(A)$ is a so-called small-world graph with a small diameter, $p_k(A)$ might be dense already for $k$ around 10.}

The approximation error that one introduces when replacing $f(A) \approx p(A)$ can conveniently be bounded using the famous result by Crouzeix and Palencia that the numerical range of a matrix is a $(1+\sqrt{2})$-spectral set~\cite{crouzeix-palencia}.
Denoting by
\[
\mathcal{W}(A) := \lbrace x^* A x :~ x\in \mathbb{C}^{n},~  \Vert x\Vert_2 = 1 \rbrace,
\]
the numerical range (or field of values) of $A$, this result implies that, when $f$ is analytic,
\begin{equation}\label{bound_crouzeix}
\Vert f(A)-p(A) \Vert_2 \leq \Vert (f-p)(A) \Vert_2 \leq \mathcal{Q} \max_{z \in \mathcal{W}(A)}\vert f(z)-p(z) \vert,
\end{equation}
for any polynomial $p$, where $\mathcal{Q} = 1$ if $A$ is normal and $\mathcal{Q} = 1+\sqrt{2}$ otherwise.\footnote{It is conjectured that $\mathcal{Q} = 2$ in the general case.} Thus, the error in the matrix function approximation $f(A) \approx p(A)$ can be bounded by the error of the polynomial approximation of the scalar function $f$ on $\mathcal{W}(A)$.

The above observations have been exploited many times in numerical algorithms for matrix functions and also have important theoretical implications: clearly, the bound~\eqref{bound_crouzeix} also holds for each individual entry of $f(A)-p(A)$. For $i \neq j$, taking \emph{any} $p \in \Pi_{\dist(i,j)-1},$ we have $[p(A)]_{ij} = 0$ so that we obtain
\begin{equation*}\label{bound}
\vert [f(A)]_{ij} \vert = \vert [f(A)-p(A)]_{ij} \vert \leq \mathcal{Q} \min_{p \in \Pi_{\dist(i,j)-1}} \max_{z \in \mathcal{W}(A)}\vert f(z)-p(z) \vert.
\end{equation*}
This means that the magnitude of the entries of $f(A)$ decays at least at the rate of the best polynomial approximation to $f$ as one moves away from the sparsity pattern of $A$~\cite{Ben_Boit, Ben_Raz, FrommerSchimmelSchweitzer2018a}.\footnote{For many functions of interest, this rate is at least exponential in the polynomial degree.} 

This fact is the basis for many efficient algorithms for approximating $f(A)$, but also for related quantities like its diagonal or trace. These methods come in two distinct flavors: on the one hand, there are methods based on truncated series representations of $f(A)$ (e.g., Faber or Chebyshev series)~\cite{Ben_Raz,Goedecker}, and on the other hand, there are probing methods~\cite{benzi2023computation, frommer2025analysis ,Frommer-Schimmel-Schweitzer,stathopoulos2013hierarchical,tang2012probing} based on graph colorings. By coloring the graph $G(A^k)$ (or equivalently, by computing a distance-$k$ coloring of $G(A)$), these methods identify which entries of $A$ play a role in computing which entries of $f(A)$, thus allowing for decoupling certain computations to improve efficiency. In many cases (unless $A$ is banded or its graph is a regular grid), the computation of $A^k$ for obtaining an appropriate coloring is the most costly step in these algorithms; see, e.g.,~\cite{benzi2023computation}. We discuss probing methods in more detail in Section~\ref{probing}.
\section{Nonzero-diagonal patterns of matrix powers}\label{sec3}
As mentioned in the previous section, the nonzero structure of powers of a banded matrix is explicitly known (assuming for simplicity that the band is fully populated), while, on the other hand, determining the sparsity structure of powers of a general sparse matrix involves certain computations that might be rather costly when $n$ is large, as they are either based on at least symbolically computing powers of $A$ or on graph traversals of $G(A)$.

\begin{figure}[t]
\centering
\begin{subfigure}{.32\textwidth}
\includegraphics[width=.95\textwidth]{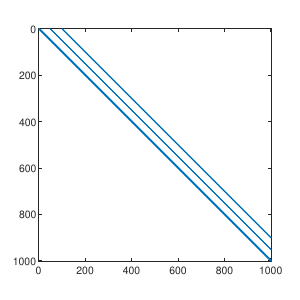}
\caption{}\label{a_fig1}
\label{Fig-A}
\end{subfigure}
\begin{subfigure}{.32\textwidth}
\includegraphics[width=.95\textwidth]{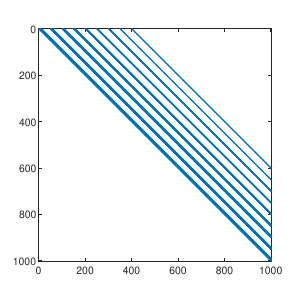}
\caption{}
\label{Fig-A4}
\end{subfigure}
\begin{subfigure}{.32\textwidth}
\includegraphics[width=.95\textwidth]{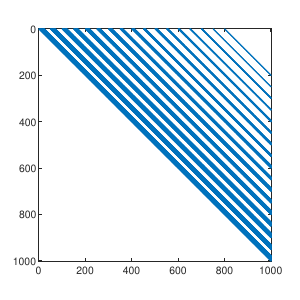}
\caption{}
\label{Fig-A8}
\end{subfigure}
\caption{Illustration of Example~\ref{example:diagonal_spread}. The matrix $A \in \mathbb{R}^{1000 \times 1000}$ contains diagonals corresponding to $\mathrm{ND}(A) = \{0,1,2,3\}\ \cup\ \{50\}\ \cup \ \{100\}$. We depict the sparsity pattern of ($\mathrm{a}$) $A$, (b) $A^4$ and (c) $A^8$.}
\label{fig:nondiag_spread}
\end{figure}

We are interested in a situation that lies between these two extremes. In certain applications, the matrix $A$ is such that only a rather small number of sub- and superdiagonals contain nonzero entries, but they are not necessarily located consecutively in the matrix (as they would be in the banded case). Examples of this situation are batch Markovian arrival process~\cite{artalejo2010markovian} in which only certain batch sizes occur or quantum walks on circulant graphs which allow for jumps of certain lengths~\cite{todtli2016continuous,yu2024classification}. In these cases, it is certainly possible to simply consider $A$ as a banded matrix considering a larger bandwidth and ignore the fact that there are ``gaps'' in the diagonal pattern, and this can lead to a severe overestimation of the actual nonzero pattern, as we demonstrate in the following example.

\begin{example}\label{example:diagonal_spread}
Let $A \in \mathbb{R}^{n \times n}$ with $\mathrm{ND}(A) = \{0,1,2,3\}\ \cup\ \{50\}\ \cup \ \{100\}$. Then $A$ can be considered as a matrix with upper bandwidth 100, which  implies that $A^4$ has an upper bandwidth of 400 and $A^8$ has an upper bandwidth of 800. This does of course not do the actual structure justice, as the band in $A$ only slowly fills up as $k$ increases and contains many ``holes''. The actual structures of the matrices (for $n = 1000$) are depicted in Figure~\ref{fig:nondiag_spread}.
\end{example}

Instead, we exploit the fact that due to the regular structure, the propagation of nonzero diagonals in $A^k$ can still be obtained very efficiently, as we will outline next. To do so, for $A \in \mathbb{C}^{n\times n}$, define the sets
\begin{equation}\label{eq:def_Sk}
\mathcal{S}_0(A) := \{ 0 \}, \qquad \mathcal{S}_k(A) := \big(\mathcal{S}_{k-1}(A)+\mathrm{ND}(A)\big) \cap \mathcal{D}_n, ~~ k\geq 1.
\end{equation}
We begin by providing a relationship between the nonzero diagonals of $A^k$ and $\mathcal{S}_k(A)$. This result is very similar in spirit to the well-known result that $A^k$ only has nonzeros in positions $(i,j)$ that are connected by length-$k$ walks in the graph of $A$. We present its short proof here for clarity of exposition.
\begin{lemma}\label{lemAB1} For $A\in \mathbb{C}^{n\times n}$,
we have $\mathrm{ND}(A^k) \subseteq \mathcal{S}_k(A)$ for all $k\geq 0$.
\begin{proof} 
We prove the result by induction. The assertion clearly holds for $k=0,1$. Now, suppose that the assertion holds for some $k$ and let $r\in \mathrm{ND}(A^{k+1})$. Then, there are $i, j \in\{1,..., n\}$ such that $r = j-i$ and $[A^{k+1}]_{ij} \neq 0$. This can be expressed as
\[ 
[A^{k+1}]_{ij} = \sum_{\ell=1}^n [A^{k}]_{i\ell}[A]_{\ell j}. 
\] 
Thus, there exists an $\ell_0$ such that $[A^{k}]_{i\ell_0}[A]_{\ell_0j} \neq 0$, which implies that both $\mathtt{diag}(A^{k},\ell_0-i)$ and $\mathtt{diag}(A,j-\ell_0)$ are nonzero. By the induction hypothesis, we have $\ell_0-i \in \mathcal{S}_{k}(A)$. Then by noting that $\mathcal{S}_{k}(A) = \left\{ i\in \mathcal{D}_n : ~i = r_1+\cdots+r_k, ~~r_j\in \mathrm{ND}(A)\right\}$, it means there are $r_1,\ldots,r_k\in \mathrm{ND}(A)$, where
$ \ell_0-i = r_1+ \cdots +r_k$. It is now evident that $j-\ell_0+ \ell_0-i\in \mathcal{S}_{k+1}(A)$, which completes the proof.
\end{proof} 
\end{lemma}

\subsection{Efficiently computing $\mathcal{S}_k(A)$}\label{FFT:Sk:new}
It depends on the size $m = |\mathrm{ND}(A)|$ which method for computing the sets $\mathcal{S}_k(A)$ is most efficient. For example, if both $m$ and $k$ are small, then it might be beneficial to compute the sets $\mathcal{S}_k(A)$ recursively from the definition~\eqref{eq:def_Sk}, requiring $\mathcal{O}(m^k)$ operations. Often, however, it is more efficient to exploit the fast Fourier transform (FFT) as follows:
Let $c_1$ denote the indicator vector of $\mathrm{ND}(A)$, i.e.,
\[
c_1(r) = \begin{cases} 
1, & \text{if } r\in\mathrm{ND}(A), \\
0, & \text{otherwise},
\end{cases}
\]
and denote by $c_k$ the $k$-fold discrete convolution of $c_1$, i.e., $c_k = c_{k-1}*c_1$, where the discrete convolution is defined by
\[
(u*w)(i)=\sum_{j\in\mathbb{Z}}u(j)w(i-j)
\]
and can efficiently be computed via the FFT.
Then the vector $c_k$ counts the number of representations of each integer as a sum of $k$ elements of $\mathrm{ND}(A)$. Consequently,
\[
\mathcal{S}_k(A)
= \{\,i\in\mathcal{D}_n:c_k(i)\neq0\,\}.
\]
If $\mathrm{ND}(A)$ contains negative integers, then the indices are first shifted to become nonnegative before applying the FFT, and the shift is removed after the computation to recover the original indices.
Using this approach, the sets $\mathcal{S}_1(A),\dots,\mathcal{S}_k(A)$ can be computed with overall complexity $\mathcal{O}(k n \log n)$. Due to the very sophisticated and well-optimized software packages being available for the FFT, this approach is often extremely fast (and in all our experiments reported in Section~\ref{sec7}, its computational complexity is completely negligible compared to the required computations in other  algorithms). 

\begin{remark}
We mention in passing that another approach for possibly even more efficiently computing the sets $\mathcal{S}_k(A)$ can be based on representing them as bit vectors and then working with bit operations performed directly in hardware. Associating each integer $\ell \in \mathcal{D}_n$ with a bit position $i = \ell + (n-1) \in \{0,1,\dots,2n-2\}$, we can define unsigned bit vectors $b_{k-1}\in\{0,1\}^{2n-1}$ whose $i$th bit is 1 if and only if $\ell = i-(n-1)\in \mathcal{S}_{k-1}(A)$. Denoting $\mathrm{ND}(A) = \{r_1,\dots,r_m\}$, the indicator vector $b_k$ of $\mathcal{S}_{k}(A)$ is then given by 
\[
\big( (b_{j-1} \ll r_1)\ |\ (b_{j-1} \ll r_2)\ |\ \dots\ |\ (b_{j-1} \ll r_m) \big)\ \mathtt{\&} \ (2^{2n-1}-1),
\]
where the notations $\mathtt{|},\mathtt{\&},\ll,$ and $\gg$ stand for bit-wise ``or'', bit-wise ``and'', left bit shift and right bit shift, respectively. As we focus on MATLAB for our implementation, and MATLAB does not support a native packed bit-vector type, we do not pursue this approach further here.
\end{remark}
\section{Identifying relevant principal submatrices of $A$}\label{sec4}
The relations given as~\eqref{eq:def_Sk} identify the indices of diagonals of $A^k$ which may contain nonzero elements, while all other diagonals are guaranteed to be zero. Clearly, if $p_k \in \Pi_k$ is a polynomial, then only the diagonals of $p_k(A)$ with indices contained in $\bigcup_{\ell=0}^k \mathcal{S}_\ell(A)$ can contain nonzero elements. 

The other way around, only certain entries of $A$ influence the computation of an individual entry $[p_k(A)]_{ij}$, while many other entries do not enter the computation at all (and therefore, can be assumed to be zero for the sake of this computation). We make this intuition precise in Lemma~\ref{submatrix_lem} below, which requires some additional notation.
Define the auxiliary sets
\begin{equation}\label{index}
\delta_{ij}^{(k)}:= \bigcup_{s=0}^{k}\bigcup_{l=0}^{s} \Big( \big(i+\mathcal{S}_l(A)\big)\cap \big(j-\mathcal{S}_{s-l}(A)\big)\Big),~~~~~\Delta_{ij}^{(k)} := \delta_{ij}^{(k)} \cap \{1,\ldots,n\}.
\end{equation}
The set $\Delta_{ij}^{(k)}$ is a superset of the set $W_{ij}^{(k)}$ of all nodes in the graph of $A$ which lie on at least one walk of length at most $k$ connecting $i$ and $j$. If all nonzero diagonals of $A$ are fully populated, we actually have $\Delta_{ij}^{(k)} = W_{ij}^{(k)}$, because the definition~\eqref{index} essentially gives the indices of all ``structural nonzeros'' of $A^k$ under this assumption. If $A$ also contains some zero entries on its nonzero diagonals, we have $W_{ij}^{(k)} \subset \Delta_{ij}^{(k)}$. Thus, $\Delta_{ij}^{(k)}$ identifies rows and columns of $A$ that may have an influence on the value of $[p_k(A)]_{ij}$. We note that either $\Delta_{ij}^{(k)} = \emptyset$ or \emph{both} $i$ and $j$ lie in $\Delta_{ij}^{(k)}$.  To see this, assume that $\Delta_{ij}^{(k)}\neq\emptyset$. Then, by definition \eqref{index}, there exists an $s$ such that $j-i\in\mathcal{S}_s(A)$. Consequently, the term corresponding to $l=0$ in the defining intersection contains $i$, while the term corresponding to $l=s$ contains $j$.
\begin{lemma}\label{submatrix_lem}
Let $A \in \mathbb{C}^{n \times n}$, $i, j \in\{1,\ldots, n\}$ and let $\delta_{ij}^{(k)}$ and  $\Delta_{ij}^{(k)}$ be defined as in~\eqref{index}. If $\Delta_{ij}^{(k)} = \emptyset$, then $[p(A)]_{ij}=0$, for any $p \in \Pi_k$; otherwise, 
\begin{equation}\label{eq:entry_polynomial_submatrix}
[p(A)]_{ij} = \left[p\left(B^{(k)}_{ij} \right)\right]_{\mathbb{O}(i)\mathbb{O}(j)},
\end{equation}
where $B^{(k)}_{ij} := A[\Delta_{ij}^{(k)}, \Delta_{ij}^{(k)}]$ and $\mathbb{O}$ is the order function on $\Delta_{ij}^{(k)}$.
\end{lemma}
\begin{proof}
We first establish the result for monomials $ A^s $, $s \geq 0 $. It is clear that if $s=0~\text{or}~1$ and $\Delta_{ij}^{(s)} \neq \emptyset$, then
$
[A^s]_{ij}
=
\bigl[(B^{(s)}_{ij})^s\bigr]_{\mathbb{O}(i)\mathbb{O}(j)}.
$
Now, let $s \ge 2$. We then have
\begin{equation}\label{pro}
\big[A^{s} \big]_{ij} = \sum_{t_1-i\in \mathrm{ND}(A)}\sum_{t_2-t_1\in \mathrm{ND}(A)}\cdots \sum_{j-t_{s-1}\in \mathrm{ND}(A)} [A]_{it_1}[A]_{t_1 t_2}\cdots [A]_{t_{s-1}j}.
\end{equation}
First, for $\ell\in\{1,\ldots,s-1\}$, we observe that for each $t_\ell$ occurring in one of the sums in \eqref{pro}, $t_\ell - i$ belongs to $ \mathcal{S}_\ell(A)$ and $j - t_\ell$ belongs to $\mathcal{S}_{s-\ell}(A)$. This implies that $t_\ell$ lies in the intersection of $i + \mathcal{S}_\ell(A)$ and $j - \mathcal{S}_{s-\ell}(A)$. Consequently, by defining the set
\[
\mathcal{G}_{ij}^{(s)}:=\bigcup_{r=0}^{s} \Big( \big(i+\mathcal{S}_{r}(A)\big)\cap \big(j-\mathcal{S}_{s-r}(A)\big)\Big)\cap \{1,\ldots,n\},
\]%
for all indices $t_\ell$, for $\ell\in\{1,\ldots,s-1\}$, we have $t_\ell\in \mathcal{G}_{ij}^{(s)}$. For $s \geq 0$, consider the principal submatrix $B^{(s)}_{ij} := A[\mathcal{G}_{ij}^{(s)},\mathcal{G}_{ij}^{(s)}]$. Therefore, concerning the relation between the entries of $B^{(s)}_{ij}$ and $A$ for $(t_{\ell}, t_{\ell+1})$, with $\ell\in\{1,\ldots,s-2\}$, for order function $\mathbb{O}$ on $\mathcal{G}_{ij}^{(s)}$, we have $[A]_{t_{\ell},t_{\ell+1}} = [B^{(s)}_{ij}]_{\mathbb{O}(t_{\ell}),\mathbb{O}(t_{\ell+1})}$. If $i$ and  $j$ are in $\mathcal{G}_{ij}^{(s)}$, then we have $[A]_{it_{\ell}}=[B^{(s)}_{ij}]_{\mathbb{O}(i)\mathbb{O}(t_{\ell})}$, and $[A]_{t_\ell j}=[B^{(s)}_{ij}]_{\mathbb{O}(t_\ell)\mathbb{O}(j)}$, implying that from \eqref{pro}, $[A^s]_{ij} = [(B^{(s)}_{ij})^s]_{\mathbb{O}(i)\mathbb{O}(j)}$. If $i$ or $j$ is not in $\mathcal{G}_{ij}^{(s)}$, then either $i$ or $j$ is not in the intersection of $i + \mathcal{S}_0(A)$ and $j - \mathcal{S}_{s}(A)$ and also is not in the intersection of $i + \mathcal{S}_s(A)$ and $j - \mathcal{S}_{0}(A)$. This means $j-i\not\in \mathcal{S}_{s}(A)$, and by Lemma \ref{lemAB1}, this implies $j-i\not\in \mathrm{ND}(A^s)$. Consequently, $[A^s]_{ij}=0$.
Now, we can extend this result by taking the union of $\mathcal{G}_{ij}^{(s)}$, for $s \in \{0, \ldots, k\}$, for any polynomial of degree $k$, to complete the proof.
\end{proof}

Typically, we are not genuinely interested in a matrix polynomial $p(A)$, but rather in some more general matrix function $f(A)$. The approach of Lemma~\ref{submatrix_lem} can also be utilized in this setting, yielding an approximation of the entry $[f(A)]_{ij}$. The approximation error can be quantified using the Crouzeix--Palencia theorem.

\begin{theorem}\label{ApproxThe}
Let $A\in \mathbb{C}^{n\times n}$ and for $ i,j \in\{1,\ldots,n\}$, let $\Delta_{ij}^{(k)} $ be defined as in \eqref{index}. If $\Delta_{ij}^{(k)} = \emptyset$, then we have
\begin{equation*}
\Big\vert\big[f(A)\big]_{ij} \Big\vert\leq 
\mathcal{Q}~ \min_{p\in \Pi_k} \max_{z \in \mathcal{W}(A)}\vert f(z)-p(z) \vert,
\end{equation*}
where $\mathcal{Q} = 1$, if $A$ is Hermitian, and $\mathcal{Q} = 1 + \sqrt{2}$, otherwise. If $\Delta_{ij}^{(k)} \neq \emptyset$, then, with $B^{(k)}_{ij} = A[\Delta_{ij}^{(k)}, \Delta_{ij}^{(k)}]$, and $\mathbb{O}$ being the order function on $\Delta_{ij}^{(k)}$, we have
\[
\Big\vert\big[f(A)\big]_{ij} - [f(B^{(k)}_{ij})\big]_{\mathbb{O}(i)\mathbb{O}(j)}\Big\vert\leq 
2\mathcal{Q}~ \min_{p\in \Pi_k} \max_{z \in \mathcal{W}(A)}\vert f(z)-p(z) \vert.
\]
\end{theorem}
\begin{proof}
The first case was already discussed in Section~\ref{sec2}. Thus, let $\Delta_{ij}^{(k)} \neq \emptyset$. Then, by Lemma~\ref{submatrix_lem}, we have $[p(A)]_{ij} = [p(B^{(k)}_{ij})]_{\mathbb{O}(i)\mathbb{O}(j)}$, for any $p \in \Pi_k$. Since $B^{(k)}_{ij}$ is a principal submatrix of $A$, we have $\mathcal{W}(B^{(k)}_{ij})\subseteq \mathcal{W}(A)$, which implies
\begin{equation}\label{ineq-f-p}
\max_{z \in \mathcal{W}(B^{(k)}_{ij})}\vert f(z)-p(z)\vert\leq \max_{z \in \mathcal{W}(A)} \vert f(z)-p(z)\vert.
\end{equation}
By utilizing $[p_k(A)]_{ij} = [p_k(B^{(k)}_{ij})]_{\mathbb{O}(i)\mathbb{O}(j)}$ from Lemma~\ref{submatrix_lem}, we arrive at
\begin{align*}
\Bigl| \bigl[f(A)\bigr]_{ij}
      - \bigl[f(B^{(k)}_{ij})\bigr]_{\mathbb{O}(i)\mathbb{O}(j)} \Bigr|
&\le
\Bigl| \bigl[f(A)\bigr]_{ij}
      - [p_k(A)]_{ij} \Bigr| \\
&+
\Bigl| \bigl[f(B^{(k)}_{ij})\bigr]_{\mathbb{O}(i)\mathbb{O}(j)}
      - [p_k(B^{(k)}_{ij})]_{\mathbb{O}(i)\mathbb{O}(j)} \Bigr| \\
&\le
2\mathcal{Q}
\max_{z\in\mathcal{W}(A)}
|f(z)-p_k(z)|,
\end{align*}
where the last inequality is obtained by first applying the bound
\[
\vert [f(D)-p(D)]_{ij}\vert
\leq
\mathcal{Q}\max_{z\in\mathcal{W}(D)}
\vert f(z)-p(z)\vert,
\]
which follows from \eqref{bound_crouzeix} and holds for any matrix
$D\in\mathbb{C}^{m\times m}$, and then using \eqref{ineq-f-p}. The proof is completed by taking the minimum over all
$p_k\in\Pi_k$.
\end{proof}

We now give an example of how Theorem~\ref{ApproxThe} can be used for approximating entries of matrix functions.

\begin{example}\label{example:diagonals}
Let $A \in \mathbb{R}^{n \times n}$ with
\[
\mathrm{ND}(A) = \{-154,\dots, -146\}\ \cup\ \{-3, \dots, 3\}\ \cup \ \{ 148, \dots, 152\}\ \cup \ \mathcal{N},
\]
for different choices of $\mathcal{N} \subseteq \mathcal{D}_{n}$. The nonzero elements of $A$ are chosen randomly in $[-1,1]$. We set $ k = 9 $ as the polynomial approximation degree in Theorem \ref{ApproxThe} and aim to approximate the entry $[\exp(A)]_{1500,1500}$. In figures~\ref{Fig-3-T32},~\ref{Fig-4-T32}, and~\ref{Fig-5-T32}, we can see the nonzero patterns of $A$ together with the submatrix (blue points) from Theorem~\ref{ApproxThe} used for approximating the element $[\exp(A)]_{1500,1500}$. The diagonals from $\mathcal{N}$ are displayed in green. 

First, we use $\mathcal{N} = \{388,\dots,392\}$ and $n = 3000$; see Figure~\ref{Fig-3-T32}. The submatrix is of size $269$ and the approximation error for $[\exp(A)]_{1500,1500}$ is $6.17 \cdot 10^{-11}$. Next, we consider  $\mathcal{N} = \{1228,\dots,1232\}$, again for $n = 3000$ (see Figure~\ref{Fig-4-T32}) and for $n = 8000$ (see Figure~\ref{Fig-5-T32}). In both cases, we obtain a submatrix of size $279$. The approximation error is $6.99 \cdot 10^{-10}$, for $n = 3000$, and $4.23 \cdot 10^{-11}$, for $n = 8000$. 

We observe that the chosen submatrix of course strongly dependent on the positions of the non-zero diagonals of $A$, and that from a sufficiently large value of $n$ onward, the size of the submatrix does not change for a particular $k$ and $(i,j)$, and the approximation quality does not deteriorate.

\begin{figure}[t]
\centering
\begin{subfigure}{.32\textwidth}
\includegraphics[width=.95\textwidth]{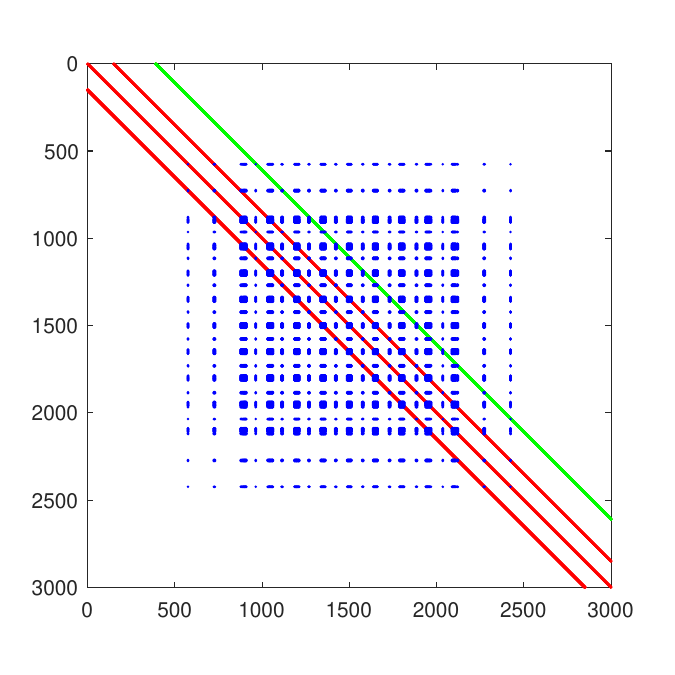}
\caption{}
\label{Fig-3-T32}
\end{subfigure}
\begin{subfigure}{.32\textwidth}
\includegraphics[width=.95\textwidth]{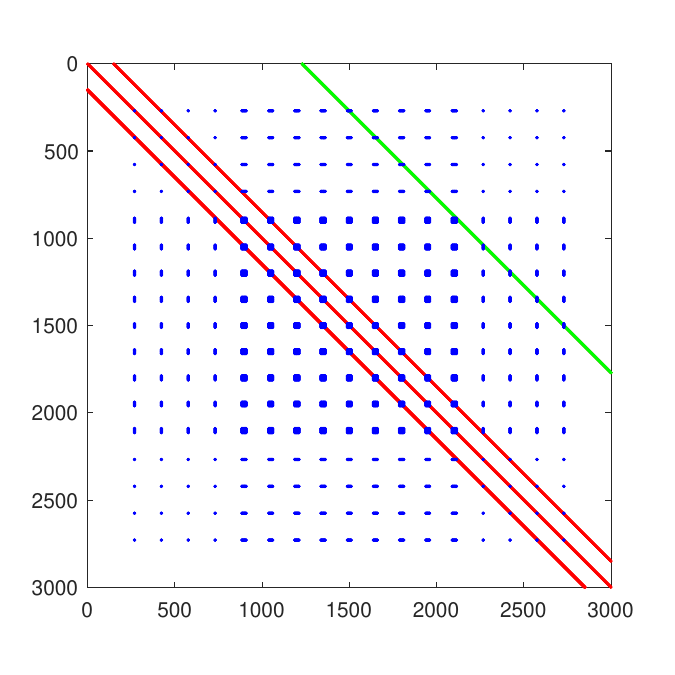}
\caption{}
\label{Fig-4-T32}
\end{subfigure}
\begin{subfigure}{.32\textwidth}
\includegraphics[width=.95\textwidth]{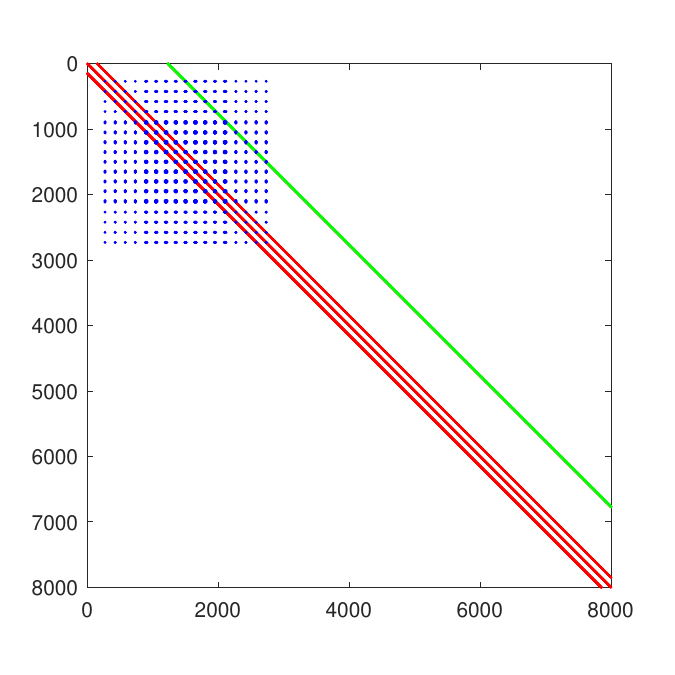}
\caption{}
\label{Fig-5-T32}
\end{subfigure}
\caption{Illustration of Example~\ref{example:diagonals}. The matrix $A$ contains diagonals corresponding to $\mathrm{ND}(A) = \{-154,\dots, -146\}\ \cup\ \{-3, \dots, 3\}\ \cup \ \{ 148, \dots, 152\}\ \cup \ \mathcal{N}$, where $\mathcal{N} = \{388,\dots,392\}$ in ($\mathrm{a}$) and $\mathcal{N} = \{1228,\dots,1232\}$ in (b) and (c). The matrix is of size $n = 3000$ in ($\mathrm{a}$) and (b) and of size $n = 8000$ in (c). The diagonals from $\mathcal{N}$ are depicted in green, the other diagonals are depicted in red. The blue dots indicate the submatrix that is used for approximating $[\exp(A)]_{1500,1500}$ using a degree-$9$ polynomial approximation.}
\end{figure}
\end{example}

\subsection{Efficient computation of $\delta_{ij}^{(k)}$ and $\Delta_{ij}^{(k)}$}\label{subsec:delta_efficient}
According to Lemma~\ref{submatrix_lem}, only the entries of a certain principal submatrix of $A$ play a role for approximating $[f(A)]_{ij}$ using a degree-$k$ polynomial approximation. Thus, if the set $\Delta_{ij}^{(k)}$ (which determines the size of the submatrix) can be efficiently computed and is small compared to the matrix size $n$, the cost for evaluating~\eqref{eq:entry_polynomial_submatrix} might be much smaller than other possible approaches. In many applications, we are not only interested in one individual entry $[p(A)]_{ij}$, but rather in many entries. We therefore now discuss how $\delta_{ij}^{(k)}$ (and thus $\Delta_{ij}^{(k)}$) can be efficiently determined for several values of $i$ and $j$ without recomputing it from scratch. To do so, we exploit several useful properties, as summarized in the following proposition.
\begin{proposition}\label{prop:facts_delta}
The sets $\delta_{ij}^{(k)}$ and  $\Delta_{ij}^{(k)}$ defined in~\eqref{index} have the following properties:
\begin{enumerate}
\item For $i \in \{1, \ldots, n\}$, we have
$
i \in \Delta_{i,i}^{(k)}.
$
\item For integer $\ell$, we have
\begin{equation}\label{T_ij+z}
\delta_{i+\ell,j+\ell}^{(k)} = \delta_{i,j}^{(k)} + \ell.
\end{equation}
\item If $\Delta_{ij}^{(k)} \subseteq M^{(k)}_{ij} \subseteq \{1, \ldots, n\}$ and $\mathbb{O}$ is the order function on $M^{(k)}_{ij}$, then we have
\begin{equation*}
[p_k(A)]_{ij} = \left[p_k\left(C^{(k)}_{ij}\right)\right]_{\mathbb{O}(i)\mathbb{O}(j)},
\end{equation*}
where $C^{(k)}_{ij} := A[M^{(k)}_{ij}, M^{(k)}_{ij}]$.

\item With 
\begin{equation}\label{eq:def_Uk}
\mathcal{U}_p(A) := \bigcup_{\ell=0}^p \mathcal{S}_{\ell}(A),    
\end{equation}
we have
\begin{equation*}
\delta_{ij}^{(k)} \subseteq \left( i + \mathcal{U}_k(A) \right) \cap \left( j - \mathcal{U}_k(A) \right).
\end{equation*}
\item Using the notation from property~4 above, the set $\delta_{ij}^{(k)}$ can equivalently be written as
\begin{equation}\label{del:uk}
\delta_{ij}^{(k)} = \bigcup_{\ell=0}^{k} \left( j - \mathcal{S}_{\ell}(A) \right) \cap \left( i + \mathcal{U}_{k - \ell}(A) \right).
\end{equation}
\end{enumerate}

\begin{proof}
We only prove part 5 here, as the first four parts are immediate.

First, let us define 
$
\mathcal{C}_{sp} := \bigl(i+\mathcal{S}_{p}(A)\bigr)\cap \bigl(j-\mathcal{S}_{s-p}(A)\bigr).
$
Now from \eqref{index} we have:
\[
\delta_{ij}^{(k)}
=
\bigcup_{0\leq p\leq s\leq k} \mathcal{C}_{sp}
=
\bigcup_{\ell=0}^{k}
\left(
\bigcup_{\substack{0\leq p,s\leq k\\s-p=\ell}}
\mathcal{C}_{sp}
\right)
=
\bigcup_{\ell=0}^{k}
\left(
\left(j-\mathcal{S}_{\ell}(A)\right)
\cap
\left(i+\mathcal{U}_{k-\ell}(A)\right)
\right),
\]
which completes the proof.
\end{proof}
\end{proposition}

In particular, part 2~of Proposition~\ref{prop:facts_delta} implies that we can compute $\delta_{i,j}^{(k)}$, for all $i, j \in \{1, \ldots, n\}$, by only computing the sets $\delta_{p,1}^{(k)}$ and $\delta_{1,p}^{(k)}$, for $p \in \{1, \ldots, n\}$, and use \eqref{T_ij+z} with $\ell \in \{0, \ldots, n-p\}$.


By using the reformulation of $\delta_{ij}^{(k)}$ in \eqref{del:uk}, having the advantage that the number of intersection and union operations is significantly reduced compared to the original definition in \eqref{index}, we can first compute $\mathcal{U}_p(A)$, for $p = 0, \ldots, k$, and then $\delta_{ij}^{(k)}$ using a recursive method for both with $\mathcal{O}(k m^k)$ operations, where $m = | \mathrm{ND}(A) |$. 
The complexity of the resulting method is $\mathcal{O}(k m^k)$ if $m^k \leq n$, and $\mathcal{O}(kn)$ otherwise.

\begin{remark}\label{rem:complexity}
The total complexity for computing all $\delta_{1,p}^{(k)}$ and $\delta_{p,1}^{(k)}$, for $p = 1, \ldots, n$, with $1 - p$ and $p - 1$ in $\mathcal{U}_k(A)$, is $\mathcal{O}(k m^k \ell)$, where $m = |\mathrm{ND}(A)|$, $\ell = |\mathcal{U}_k(A)|$.
\end{remark}

\begin{remark}\label{rem:krylov:more:one}
If only a single entry $[f(A)]_{ij}$ is to be computed, a straightforward way is to use a Krylov subspace method for approximating the bilinear form $e_i^*f(A)e_j$; see, e.g.,~\cite{golub2009matrices}. Doing $k$ iterations of such a method can be expected to yield a comparable error to the approach outlined above using the same polynomial degree $k$. The main cost of one Krylov iteration is a matrix vector product with $A$, which can be performed in $\mathcal{O}(mn)$ operations, where $m = |\mathrm{ND}(A)|$. Thus, the overall cost of this approach is $\mathcal{O}(knm)$. In light of Remark~\ref{rem:complexity}, the computation of all $\delta$-values is thus cheaper only if $m^{k-1}\ell \leq n$, which will typically only be the case for very small values of $m$ and $k$. Therefore, for a single entry, our approach will typically not be competitive. It is rather beneficial if several entries are required (as the $\delta$-sets only need to be computed once upfront, while a new call to a Krylov method is required for each entry). The main motivations for using our methodology are therefore the algorithms and applications discussed in sections~\ref{sec5} and \ref{sec6}.
\end{remark}
\section{Relation to probing methods}\label{probing}\label{sec5}
A class of methods similar in principle to the methods we discussed here are the so-called \emph{probing methods}. These methods were originally introduced for approximating the whole diagonal of a matrix function~\cite{tang2012probing}, but they can also be used for estimating its trace~\cite{benzi2023computation,stathopoulos2013hierarchical,Frommer-Schimmel-Schweitzer,frommer2025analysis}, and even for computing a sparse approximation of $f(A)$ itself~\cite{Frommer-Schimmel-Schweitzer,schimmel2019bounds}. We focus on the latter two applications here.

To form a sparse approximation satisfying an error bound similar to that of Theorem~\ref{ApproxThe} by a probing method, one first computes a distance-$2k$ coloring of the graph of $A$.\footnote{A distance-$d$ coloring of a graph is a coloring in which no two nodes $i$ and $j$ that have a geodesic distance $\text{dist}(i,j) \leq d$ are allowed to have the same color. Thus, the classical graph coloring problem corresponds to the choice $d=1$.} Such a coloring induces a partition of the node set as follows:
\begin{equation}\label{eq:node_partitioning}
\{1,\dots,n\}=P_1 \cup \ldots \cup P_m, \enspace P_\ell \neq \emptyset \mbox{ for } \ell = 1,\ldots,m  \text{ and } P_\ell \cap P_j = \emptyset \text{ for } \ell \neq j,
\end{equation}
where $P_\ell$ is the set of all nodes with assigned color $\ell$. The \emph{probing vectors} corresponding to the partitioning~\eqref{eq:node_partitioning} are then defined as
\begin{equation}\label{eq:probing_vectors}
v_\ell:=\sum\limits_{i\in P_\ell} e_i, \; \ell \in \{1,\ldots,m\},
\end{equation}
where $e_i$ is the $i$th canonical unit vector. Based on these probing vectors, a sparse approximation $\widetilde{F}$ of $f(A)$ can be computed by
\begin{equation*}
 \widetilde{F}_{ij}:=
 \begin{cases}
    [f(A)v_\ell]_i \text{ for } j \in V_\ell, &\text{ if } \text{dist}(i,j) \leq k, \\
    0, &\text{ if } \text{dist}(i,j) > k.
    \end{cases}
\end{equation*}
Given the partitioning~\eqref{eq:node_partitioning}, forming the sparse approximation $\widetilde{F}$ thus requires computing $m$ matrix vector products with $f(A)$, which can, for example, be approximated by a Krylov subspace method, and then distributing the entries of the resulting vectors according to~\eqref{eq:probing_vectors}. For details on the derivation of this method, see~\cite{Frommer-Schimmel-Schweitzer,schimmel2019bounds}. 

Except for certain highly structured special cases (e.g., banded matrices or matrices with an underlying regular grid structure), the most costly part of a probing method is the first phase, in which the distance-$2k$ coloring of the graph of $A$ is computed; see also~\cite{benzi2023computation} for a discussion of this topic. Typically, the distance-$2k$ coloring is computed by forming the matrix $A^{2k}$ and then using a greedy algorithm to compute a classical (i.e., distance-$1$) coloring of its graph. The precise cost of this approach highly depends on the sparsity pattern of $A$, but it is often prohibitive even for moderate values of $k$.

Our approach outlined in sections~\ref{sec2}--\ref{sec4} can be seen as an alternative for inferring information about the sparsity structure of powers of $A$ while avoiding the explicit computation of $A^{2k}$, as we discuss next.


\subsection{Computing a node partitioning based on our methodology}\label{sec51}

Algorithm~\ref{alg:partition} gives one possible way of obtaining a partition $\mathcal{P} = \{ P_1,\ldots,P_m\}$ based on our methodology developed in sections~\ref{sec3} and \ref{sec4} with a computational complexity of $\mathcal{O}(n)$.

\begin{algorithm}
\caption{Construction of the partition $\mathcal{P} = \{P_1, \ldots, P_m\}$}\label{alg:partition}
\begin{algorithmic}[1]
\STATE Set $N = \{1, \ldots, n\}$ and $\ell = 0$.
\WHILE{$N \neq \emptyset$}
    \STATE $\ell \gets \ell + 1$.
    \STATE Initialize $P_{\ell} = \emptyset$.
    \STATE Set $G \gets N$.
    \WHILE{$G \neq \emptyset$}
        \STATE Choose $i \in G$ .
        \STATE Add $i$ to $P_{\ell}$.
        \STATE Remove $i$ and all elements of $i + \mathcal{U}_k(A)$ from $G$.
    \ENDWHILE
    \STATE Remove all elements of $P_{\ell}$ from $N$.
\ENDWHILE
\STATE \textbf{return} $\mathcal{P} = \{P_1, \ldots, P_m\}$.
\end{algorithmic}
\end{algorithm}

\begin{remark}\label{rem:pAij0}
By construction, the partition $\mathcal{P} = \{P_1, \ldots, P_m\}$ in Algorithm~\ref{alg:partition} is such that for any distinct $i, j \in P_{\ell}$, we have $\Delta_{ij}^{(k)} = \emptyset$. This is due to the fact that in line~9 of Algorithm~\ref{alg:partition}, all elements of $i + \mathcal{U}_k(A)$ are removed from the candidate set $G$. Thus, all $j \neq i$ in $P_\ell$ must satisfy $j-i\not\in \mathcal{U}_k(A)$. Consequently, for any polynomial $p_k$ of degree $k$, it follows that $[p_k(A)]_{ij} = 0$.
\end{remark}

\subsection{Delta trace estimator}\label{sec52}
Given the partition $\mathcal{P} = \{P_1, \ldots, P_m\}$ obtained from Algorithm~\ref{alg:partition}, let $S_\ell \subseteq \{1,\dots,n\}, \ell = 1,\dots,m,$ be sets such that
\begin{equation*}
\Delta_{P_{\ell}} = \bigcup_{i\in P_{\ell}} \Delta_{ii}^{(k)} \subseteq S_{\ell},
\end{equation*}
and define the corresponding \emph{delta vectors},
\begin{equation}\label{vec:prob:new}
w_{\ell}  = \sum_{i\in P_{\ell}} e_{\mathbb{O}_{\ell}(i)},~~~~~~~~\ell \in \{1,\ldots,m\},
\end{equation}
where $\mathbb{O}_{\ell}$ is the order function on $S_{\ell}$. If we simply choose $S_\ell = \{1,\dots,n\}$, for all $\ell$, then the vectors $w_\ell$ from~\eqref{vec:prob:new} exactly coincide with the probing vectors from~\eqref{eq:probing_vectors}. We are free to choose smaller sets $S_\ell$,  with the other natural choice being $S_\ell = \Delta_{P_\ell}$, for all $\ell$. In the following, we denote the submatrices of $A$ corresponding to our choice of $S_\ell$ by
\begin{equation}\label{BP:submat:Del}
B_{P_{\ell}} = A[S_{\ell},S_{\ell}].
\end{equation}
Note that for the trivial choice $S_\ell = \{1,\dots,n\}$, we have $B_{P_\ell} = A$.

An approximation of the trace of $f(A)$ now arises via the \textit{delta-trace approximation},
\begin{equation}\label{delta:trace:approx}
\mathcal{T}_{\Delta}^{(k)}(f(A)) := \sum_{\ell=1}^m w_{\ell}^* f(B_{P_{\ell}}) w_{\ell}.
\end{equation}
For future reference, we refer to the approximation~\eqref{delta:trace:approx} using  $S_\ell = \{1,\dots,n\}$ as \emph{``full''} method (this is essentially the \textit{probing-trace approximation} extensively studied in \cite{benzi2023computation,stathopoulos2013hierarchical,Frommer-Schimmel-Schweitzer,frommer2025analysis}, but with a different choice of the partition $\mathcal{P}$), and to the approximation~\eqref{delta:trace:approx} using $S_{\ell} = \Delta_{P_{\ell}}$ as \emph{``split''} method. 

The error of the approximation~\eqref{delta:trace:approx} is analyzed using the following lemma and theorem.
\begin{lemma}\label{lem:vfAv}
Let $A\in \mathbb{R}^{n\times n}$ and $B_{P_{\ell}}$ be the submatrix of $A$ as defined in \eqref{BP:submat:Del}. Let $w_{\ell}$ be the delta vector defined in \eqref{vec:prob:new} corresponding to the set $P_{\ell}$ in the partition $\mathcal{P}$ from Algorithm~\ref{alg:partition}. Then, we have
 \[
\Big\vert w_{\ell}^* f(B_{P_{\ell}}) w_{\ell} - \sum_{i\in P_{\ell}} [f(A)]_{ii}  \Big\vert \leq  2\mathcal{Q} |P_{\ell}| \min_{p_k\in \Pi_k} \max_{z \in \mathcal{W}(A)}\vert f(z)-p(z) \vert,
 \]
where $\mathcal{Q} = 1$, if $A$ is Hermitian, and $\mathcal{Q} = 1 + \sqrt{2}$, otherwise.
\begin{proof}
First, it is clear that for any polynomial $p_k$, we have 
$\sum_{i\in P_{\ell}} [p_k(A)]_{ii} = w_{\ell}^* p_k(B_{P_{\ell}}) w_{\ell}$, because for distinct $i,j\in P_{\ell}$, $[p_k(A)]_{ij}=0$ (cf.~Remark~\ref{rem:pAij0}). Thus, we have
\[
\Big\vert w_{\ell}^* f(B_{P_{\ell}}) w_{\ell} - \sum_{i\in P_{\ell}} [f(A)]_{ii} \Big\vert
\leq 
\Big\vert \sum_{i\in P_{\ell}} [f(A)-p_k(A)]_{ii} \Big\vert
+ 
\Big\vert w_{\ell}^* \left(f(B_{P_{\ell}}) - p_k(B_{P_{\ell}})\right) w_{\ell} \Big\vert.
\]
Analogous to the proof of Theorem~\ref{ApproxThe}, we can arrive at the bounds,
\[
\vert \sum_{i\in P_{\ell}} [f(A)-p_k(A)]_{ii} \vert
\leq \mathcal{Q}|P_{\ell}| \max_{z \in \mathcal{W}(A)}\vert f(z)-p(z) \vert
\]
and 
\[\vert w_{\ell}^* \left(f(B_{P_{\ell}}) - p_k(B_{P_{\ell}})\right) w_{\ell} \vert
\leq 
\mathcal{Q}\Vert  w_{\ell} \Vert _2^2 \max_{z \in \mathcal{W}(A)}\vert f(z)-p(z) \vert, 
\]
from which the desired result follows immediately.
\end{proof}
\end{lemma}
By applying Lemma~\ref{lem:vfAv} to all the sets $B_\ell, \ell = 1,\dots, m$, we directly obtain an error bound for the trace approximation~\eqref{delta:trace:approx}. This is similar to error bounds available for probing methods (see, e.g.,~\cite{Frommer-Schimmel-Schweitzer,schimmel2019bounds}), but holds in the more general setting that the $B_\ell$ do not necessarily need to agree with $A$.
\begin{theorem}\label{theorem:err:trace}
Under the assumptions of Lemma~\ref{lem:vfAv}, we have
\[
\vert \mathrm{trace}(f(A)) - \mathcal{T}_{\Delta}^{(k)}(f(A)) \vert \leq 2\mathcal{Q}n \min_{p_k\in \Pi_k} \max_{z \in \mathcal{W}(A)}\vert f(z)-p(z) \vert,
\]
where $\mathcal{Q} = 1$, if $A$ is Hermitian, and $\mathcal{Q} = 1 + \sqrt{2}$, otherwise.
\begin{proof}
The result directly follows by applying Lemma~\ref{lem:vfAv} for $\ell = 1, \ldots, m$.
\end{proof}
\end{theorem}
\begin{remark}
By Theorem~\ref{theorem:err:trace}, the same error bound is valid for both the ``full'' and the ``split'' probing methods. In practice, one can observe that the actual error of the ``split'' method is typically slightly larger than that of the ``full'' method, due to the additional error introduced the use of the decomposition into the submatrices $B_{P_{\ell}}$. The ``split'' method can often be more cost-effective, though, since it requires less storage and works with smaller matrices.
\end{remark}

\subsection{Stochastic delta estimator} 
As analyzed in \cite{frommer2025analysis}, probing methods can be combined with the stochastic Hutchinson estimator~\cite{Hutch} (instead of using the ``deterministic'' probing vectors~\eqref{eq:probing_vectors}) to improve efficiency. Of course, our approach based on delta sets can easily be combined with this \textit{stochastic probing estimator}. We describe this approach for symmetric $A$ (which is the typical situation for trace estimation). To simplify notation, in the following we let $[M]_{\mathcal{I}}$ denote the submatrix of $M$ constructed from the rows and columns with indices in $\mathcal{I}$.

Using the notations introduced in section~\ref{sec51} and \ref{sec52}, we define the \textit{stochastic delta vectors},
\[
w_{\ell}^{(p)} = \sum_{i \in P_{\ell}} X_i e_{i'},~~~~~~p = 1,\ldots,N_{\ell}
\]
where $N_\ell \geq 1$ is the number of \emph{samples} and the $X_i$ are i.i.d. Rademacher random variables such that $\mathbb{E}[X_i] = 0$ and $\mathbb{E}[(X_i)^2] = 1$. Based on these stochastic delta vectors, we introduce the \textit{stochastic delta trace} estimator as follows:
\begin{equation}\label{del:Stoch:Estim}
\mathcal{RT}_{\Delta}^{(k)}(f(A)) := \sum_{\ell=1}^m\dfrac{1}{N_{\ell}}\sum_{p=1}^{N_{\ell}}\big(w^{(p)}_{\ell}\big)^*f(B_{P_{\ell}}) w^{(p)}_{\ell}.
\end{equation}
The expected value of this estimator is
\begin{equation*}
\mathbb{E}\left[\mathcal{RT}_{\Delta}^{(k)}(f(A))\right] = \sum_{\ell = 1}^m \sum_{i\in P_{\ell}} [f(B_{P_{\ell}})]_{\mathbb{O}_{\ell}(i)\mathbb{O}_{\ell}(i)}.
\end{equation*}
Thus, for the ``full'' method, the expected value exactly equals $\mathrm{trace}(f(A))$, while for the ``split'' method, it is only approximately equal to $\mathrm{trace}(f(A))$, with an approximation error as in Theorem \ref{theorem:err:trace}.\footnote{In practice, the error introduced by the stochastic estimator is typically larger than the approximation error resulting from the decomposition of $A$ into submatrices, which is therefore typically negligible.} As before, the ``full'' method coincides with the standard (stochastic) probing method of~\cite[Section 3]{frommer2025analysis}. As in Theorem~\ref{theorem:err:trace} for the error of the deterministic method, we can again straightforwardly generalize the analysis of the method to the ``split'' version. In particular, both for the full and split methods, for every $\varepsilon > 0$ we have the tail bound as
\[
\mathbb{P}\left( \left|\operatorname{trace}(f(A)) - \mathcal{RT}_{\Delta}^{(k)}(f(A)) \right| \geq \varepsilon \right) \leq 2 \exp\left( - \frac{\varepsilon^2}{8 \eta_1 + 8 \varepsilon \eta_2}\right),
\]
where 
\[
\eta_1 := \sum_{\ell=1}^{m} \frac{1}{N_{\ell}} \Vert M_{\ell} \Vert_F^2 ,~~~~~~~~~~~~~
\eta_2 := \max_{\ell=1,\ldots,m} \frac{1}{N_{\ell}} \|  M_{\ell} \|_2,
\]
with 
\[
M_{\ell}:=[f(B_{P_{\ell}})]_{P_{\ell}} - 
\begin{bmatrix}
\left[\mathtt{diag}\left([f(A)]_{P_{\ell}}\right)\right]_1 & & 0\\
& \ddots &\\
0 & & \left[\mathtt{diag}\left([f(A)]_{P_{\ell}}\right)\right]_{|P_{\ell}|}
\end{bmatrix};
\]
cf.~\cite[Theorem~3.8]{frommer2025analysis}.

Furthermore, as in~\cite[Lemma~3.5]{frommer2025analysis}, the \textit{variance} of the estimator \eqref{del:Stoch:Estim} is given by
\begin{equation}\label{varia:VL}
\mathbb{V}\!\left[\mathcal{RT}_{\Delta}^{(k)}(f(A))\right]
= \sum_{\ell=1}^m \frac{V_{\ell}}{N_{\ell}},
\end{equation}
where
\begin{equation}\label{varia:VL2}
V_{\ell} := \mathbb{V}\!\left[\big(w_{\ell}^{(p)}\big)^* f(B_{P_{\ell}})\, w_{\ell}^{(p)}\right]
= 2 \sum_{\substack{i,j \in P_{\ell} \\ i \neq j}} \big|[f(B_{P_{\ell}})]_{\mathbb{O}_{\ell}(i)\mathbb{O}_{\ell}(j)}\big|^2.
\end{equation}

\begin{remark}
It can be seen from~\eqref{varia:VL} and \eqref{varia:VL2} (and it is intuitively clear) that in order to achieve a low overall variance of the trace estimator, one should choose a larger number of samples $N_\ell$ when the variance of the estimator for the submatrix $B_\ell$ is large, while a small number of samples suffices if the variance of the estimator for $B_\ell$ is small. A suitable number of samples can thus be chosen either based on a priori bounds for the variances (see also Theorem~\ref{thm:new_bound_variance} below) or by estimating the variances on the fly during the method, based on computed quantities; see~\cite[sections~3.2 and~3.4]{frommer2025analysis}.
\end{remark}

The following result gives an a priori bound for the variance of our estimator. This improves upon the corresponding result of~\cite[Proposition~3.10]{frommer2025analysis} in two ways: first, it is slightly sharper (the constant in the bound is reduced from~$4$ to $2$) and second, it holds under more general conditions (in~\cite[Proposition~3.10]{frommer2025analysis}, an assumption on the sign of the entries of $f(A)$ is necessary).
\begin{theorem}\label{thm:new_bound_variance}
Let $A \in \mathbb{R}^{n\times n}$ be symmetric, and let $\mathcal{P} = \{P_1, \ldots, P_m \}$ be the partition obtained from Algorithm~\ref{alg:partition}. Then, $V_{\ell}$, as defined in \eqref{varia:VL}, satisfies
\[
V_{\ell} \leq 2 |P_{\ell}| \min_{p_k\in\Pi_k} \max_{z \in \mathcal{W}(A)}\vert f(z)-p(z) \vert^2.
\]
\begin{proof}
First, for $p_k\in \Pi_k$ and $i,j\in P_{\ell}$,  $i\neq j$, we have $[p_k(B_{P_{\ell}})]_{\mathbb{O}_{\ell}(i)\mathbb{O}_{\ell}(j)} = 0$. Therefore, from~\eqref{varia:VL2} and since the main diagonal entries are not included in the sum $V_\ell$, we have
\begin{align*}
\dfrac{V_{\ell}}{2} &= \sum_{\substack{i,j \in P_{\ell} \\ i \neq j}} \big|[f(B_{P_{\ell}})-p_k(B_{P_{\ell}})]_{\mathbb{O}_{\ell}(i)\mathbb{O}_{\ell}(j)}\big|^2 \leq \Vert [f(B_{P_{\ell}})-p_k(B_{P_{\ell}})]_{\mathbb{O}_{\ell}(P_{\ell})} \Vert_F^2 \\
&\leq |P_{\ell}|\Vert [f(B_{P_{\ell}})-p_k(B_{P_{\ell}})]_{\mathbb{O}_{\ell}(P_{\ell})} \Vert_2^2,
\end{align*}
where we used the fact that for an $m \times m$ matrix $C$, we have $\|C\|_F^2 \leq m\|C\|_2^2$. Furthermore, as $[f(B_{P_{\ell}})-p_k(B_{P_{\ell}})]_{\mathbb{O}_{\ell}(P_{\ell})}$ is a principal submatrix of $f(B_{P_{\ell}})-p_k(B_{P_{\ell}})$, it has a smaller spectral norm, so that
\[
|P_{\ell}|\Vert [f(B_{P_{\ell}})-p_k(B_{P_{\ell}})]_{\mathbb{O}_{\ell}(P_{\ell})} \Vert_2^2 \leq |P_{\ell}|\Vert f(B_{P_{\ell}})-p_k(B_{P_{\ell}}) \Vert_2^2.
\]
Since $ \Vert f(B_{P_{\ell}})-p_k(B_{P_{\ell}}) \Vert_2^2 \leq \max_{z \in \mathcal{W}(A)}\vert f(z)-p(z) \vert^2$, the result follows by taking the minimum over all $p_k\in \Pi_k$.
\end{proof}
\end{theorem}
\begin{remark}
The stochastic delta estimator defined by \eqref{del:Stoch:Estim} can typically be expected to achieve a lower error than that of the deterministic delta trace approximation method at a similar computational cost. In certain situations, one can rigorously guarantee that this is the case. For example, even when only one stochastic sample is used per $P_\ell$,the stochastic error will always be below the deterministic error if $f$ and $A$ are such that $[f(A)]_{ij}$ has a constant sign for every $\ell$ and for all $i, j \in P_{\ell}$ (see \cite[Proposition 3.15]{frommer2025analysis}). 
\end{remark}
\section{An efficient algorithm for computing functions of sparse Toeplitz matrices}\label{sec6}
Here, we illustrate how our methodology can be used for computing functions of sparse Toeplitz matrices by an approach that is typically much more efficient than the probing-based approximation scheme discussed in Section~\ref{sec5}.

A Toeplitz matrix $T \in \mathbb{C}^{n\times n}$ has the form
\begin{equation*}
T = \begin{pmatrix} a_0 & a_{1} & a_{2} & \ldots & a_{(n-1)} \\
a_{-1} & a_0 & a_{1} & \ldots & a_{(n-2)} \\
a_{-2} & a_{-1} & a_0 & \ldots & a_{(n-3)} \\
\vdots & \vdots & \vdots & \ddots & \vdots \\
a_{-(n-1)} & a_{-(n-2)} & a_{-(n-3)} & \ldots & a_0 \end{pmatrix},
\end{equation*}
meaning that all entries along the diagonal $\mathtt{diag}(T,r)$ have the same value $a_r$.

It is important to note that the powers $T^k, k > 1$, of a Toeplitz matrix are in general \emph{not} Toeplitz any longer---while this property holds for (bi-)infinite Toeplitz matrices. This is due to effects incurred by the finite summation, which leads to differing entries near the ``boundaries'' of the matrix, while large parts of the ``interior'' diagonal values will still be constant, at least if only a small, fixed number of diagonals of $T$ contains nonzeros---precisely the setting of our paper. We refer to~\cite[Sections~4.3 and~5.2]{Gray2006} for a detailed discussion of this phenomenon and quantitative characterizations. Due to this effect, if $k \ll n$, \emph{many} entries along each diagonal will be equal, which can be exploited for obtaining a very efficient algorithm, as we will outline in the remainder of this section: By identifying the positions in $T^k$ that satisfy
\begin{equation}\label{cond:Tk:Toeplitz}
\left[T^k\right]_{ij} = \left[T^k\right]_{i-1,j-1},
\end{equation}
we aim to reduce the complexity of computing $f(T)$ by reducing it to the computation of a function of  a typically much smaller submatrix. This can be achieved
by investigating the structure of $\delta_{ij}^{(k)}$ for the matrix $T^k$ and identifying all possible elements that satisfy condition \eqref{cond:Tk:Toeplitz}. 

Based on this observation, it is sufficient to select only one \emph{representative entry} among all entries satisfying condition \eqref{cond:Tk:Toeplitz} and compute its value by the methods of Section~\ref{sec4}. Then, this value can be distributed to all entries represented by it in the resulting matrix.
This means that computing a polynomial in $T$ can typically be reduced to the computation of a polynomial of a single, much smaller submatrix of $T$ whose entries are then repeatedly placed in a matrix of size $n \times n$.

If two neighboring entries on the same diagonal of the matrix generate identical submatrices from their corresponding $\delta$-sets using Lemma~\ref{submatrix_lem}, and if the positions of these entries within their respective submatrices are also identical, then it is easy to observe that the corresponding entries in the function of the original matrix will also be (approximately) equal. Based on this fact, we present a condition in the next lemma.

\begin{lemma}\label{lem:toeplitz}

Assume that the condition
\begin{equation}\label{cond1}
\Delta_{i,j}^{(k)} = \Delta_{i-1,j-1}^{(k)} + 1,
\end{equation}
holds.
Then, for $\ell = 0,\dots,k$, we have
\begin{equation*}
\left[T^{\ell}\right]_{ij} = \left[T^{\ell}\right]_{i-1,j-1}.
\end{equation*}
\end{lemma}
\begin{proof}
First, we define two submatrices,
\[
T_1 := T[\Delta_{i,j}^{(k)}, \Delta_{i,j}^{(k)}],
~~~~~~\text{and}~~~~~~
T_2 := T[\Delta_{i-1,j-1}^{(k)}, \Delta_{i-1,j-1}^{(k)}].
\] 
The goal is to show that these two submatrices are equal. Let $\mathbb{O}_*$ denote the order function on $\Delta_{i-1,j-1}^{(k)}$, and let $\mathbb{O}$ denote the order function on $\Delta_{i,j}^{(k)}$.
Therefore, for every $p \in \Delta_{i,j}^{(k)}$, we have $\mathbb{O}(p) = \mathbb{O}_*(p-1).$
Moreover, for each entry in $T_1$ and its corresponding index in $T_2$, there exist $x,y \in \{1,\ldots,n\}$ such that $[T_1]_{\mathbb{O}(x),\mathbb{O}(y)} = [T]_{x,y}$ and $[T_2]_{\mathbb{O}_*(x-1),\mathbb{O}_*(y-1)} = [T]_{x-1,y-1}.$
Since $T$ is a Toeplitz matrix, we have $[T]_{x-1,y-1} = [T]_{x,y},$ which implies that $[T_1]_{\mathbb{O}(x),\mathbb{O}(y)} = [T_2]_{\mathbb{O}_*(x-1),\mathbb{O}_*(y-1)}.$ Therefore, $T_1 = T_2,$ and by applying Lemma~\ref{submatrix_lem}, the proof is complete.
\end{proof}

Note that Lemma \ref{lem:toeplitz} may involve a class of entries that one after another satisfy condition \eqref{cond1}. This means, for $q>0$ and a pair $(i,j)$ we may have:
\[
\Delta_{i,j}^{(k)} = \Delta_{i-p,j-p}^{(k)} + p, ~~~~ for~~ p = 0,\ldots,q.
\]
In this case, we only need to choose one representative entry of this class, e.g., $(i,j)$. 

According to the above discussion, for the Toeplitz matrix $T$, it is sufficient to consider only those sets $\Delta_{ij}^{(k)}$ for which, in $T^{\ell}$ with $\ell = 0,\dots,k$, the corresponding entry \textit{does not} have the same value as its neighboring entry on the same diagonal, meaning that the condition \eqref{cond1} \textit{does not} hold. In other words, our goal is to identify the set
\begin{equation}\label{Glob:del}
\Delta_G^{(k)} = \bigcup_{(i,j)\in \mathcal{P}} \Delta_{ij}^{(k)},
\end{equation}
where 
\begin{equation}\label{P:4:Glob:del}
\mathcal{P} = \left\{(i,j):  \Delta_{i,j}^{(k)} \neq \Delta_{i-1,j-1}^{(k)} + 1 \right\}.
\end{equation}

To compute $\Delta_G^{(k)}$, we can traverse each nonzero diagonal whose index belongs to $\mathcal{U}_k(T) = \bigcup_{\ell = 0}^k \mathcal{S}_{\ell}(T)$ and verify whether the condition on $\mathcal{P}$, hold. If so, we may add the corresponding delta-set to $\Delta_G^{(k)}$. Then, using $\Delta_G^{(k)}$, we can compute the function of the corresponding submatrix in the original Toeplitz matrix $T$ as
\begin{equation}\label{eq:Xf}
X_f = f\left(T\left(\Delta^{(k)}_G,\Delta^{(k)}_G\right)\right).
\end{equation}
Eventually, to obtain an approximation of $f(T) $, we can appropriately place the entries of $X_f$ into a sparse matrix of the same size as $ T$, using the ordering function $\mathbb{O}(\cdot)$ over $\Delta^{(k)}_G$ and the indices corresponding to $\mathcal{P}$, and then extend this placement to neighboring entries on the same diagonal that are not in $\mathcal{P}$. The main idea of this method is illustrated in Figure~\ref{Fig:Toep:method}.

\begin{algorithm}
\caption{Approximation of function of Toeplitz matrix}\label{alg:repeat}
\textbf{Input:} Toeplitz matrix $T\in \mathbb{C}^{n\times n}$, function $f$, polynomial degree $k$.\\
\textbf{Output:} Approximation $\widetilde{F}$ for $f(T)$.
\begin{algorithmic}[1]
\STATE Compute $\mathcal{S}_0(T),\ldots,\mathcal{S}_k(T)$ and set $\mathcal{U}_k(T)$.
\STATE Initialize $\Delta_G^{(k)} = \emptyset$.
\FOR{$r \in \mathcal{U}_k(T)$}
    \STATE Compute starting indices $(i,j)$ on the $r$-th diagonal.
    \STATE $\delta_0 \gets \delta_{ij}^{(k)}$.
    \IF{$i,j \in \Delta_{ij}^{(k)}$}
        \STATE Set $\Delta_G^{(k)} \leftarrow \Delta_G^{(k)} \cup \Delta_{ij}^{(k)}$ and mark entry $(i,j)$ as ``non-repeated''
    \ENDIF
    \FOR{$\ell = 1$ to $n - |r|$}
        \STATE $\Delta_\ell \gets (\delta_0 + \ell) \cap \{1, \ldots, n\}$.
        \IF{$(i+\ell,j+\ell)$ satisfies conditions \eqref{cond1}}
            \STATE Mark entry $(i+\ell,j+\ell)$ as ``repeated''
        \ELSE
            \STATE Set $\Delta_G^{(k)} \leftarrow \Delta_G^{(k)} \cup \Delta_\ell$, and mark entry $(i+\ell,j+\ell)$ as ``non-repeated''
        \ENDIF
    \ENDFOR
\ENDFOR
\STATE Compute $X_f = f(T(\Delta^{(k)}_G, \Delta^{(k)}_G))$.
\STATE Initialize zero matrix $\widetilde{F}$ and fill its entries marked as ``non-repeated'' from $X_f$.
\STATE Along each diagonal, copy the last ``non-repeated'' value to all the following ``repeated'' entries.
\end{algorithmic}
\end{algorithm}

\begin{figure}
\begin{center}
\begin{tikzpicture}[
>=Latex,
every node/.style={font=\small}
]

\node (A) {
\begin{tikzpicture}[scale=0.8]
\begin{scope}[rotate=50]

\fill[gray!25, rounded corners]
(0.75,-0.3) rectangle (1.25,1.5);

\foreach \i in {0,1,3,4}
\fill (1,{-0.11+0.35*\i}) circle (2pt);

\node at (1.09,0.79) {$\ddots$};
\node at (1.07,-0.6) {$\ddots$};
\node at (1.07,1.95) {$\ddots$};

\end{scope}
\end{tikzpicture}
};

\node[right=1.9cm of A] (C)
{
\footnotesize
$
\left[
f\!\left(
\begin{tikzpicture}[baseline=(DG.base),scale=0.8]
\node[
fill=gray!25,
rounded corners=2pt,
minimum width=1.5cm,
minimum height=1.8cm,
inner sep=4pt,
align=center
](DG)
{
\begin{tabular}{c}
\footnotesize submatrix of\\
$\Delta_G$
\end{tabular}
};
\end{tikzpicture}
\right)
\right]_{\mathbb{O}(\bullet)}
=
\bigstar
$
};


\node[right=1.5cm of C] (D)
{\footnotesize
\begin{tikzpicture}[scale=0.8]
\begin{scope}[rotate=50]

\fill[gray!25, rounded corners]
(0.75,-0.37) rectangle (1.25,1.55);

\foreach \i in {0,1,3,4}
{
\node at (0.825,{0.12+0.35*\i}) {\scriptsize $\bigstar$};
}
\node at (0.87,0.85) {$\ddots$};
\node at (0.85,-0.3) {$\ddots$};
\node at (0.85,2.2) {$\ddots$};

\end{scope}
\end{tikzpicture}};

\draw[->, thick] 
(C.east) -- (D.west);

\draw[->, thick]
(A.east) -- (C.west);

\node[below=5pt of A] {\scriptsize \bf ($a$)};
\node[below=5pt of C] {\scriptsize \bf (b)};
\node[below=5pt of D] {\scriptsize \bf (c)};

\end{tikzpicture}
\end{center}
\caption{
(a) The points marked by $\bullet$ represent the entries along $\mathtt{diag}(T,r)$ that satisfy condition \eqref{cond1} one after another. 
(b) First, a $\delta$-set is extracted for the corresponding selected representative points $\bullet$ of each class, then all delta sets are added to $\Delta_G$ and the matrix function of the corresponding submatrix of $\Delta_G$ is computed. 
(c) The resulting value $\bigstar$ is placed into the resulting matrix at exactly the same positions as the original points~$\bullet$.
}\label{Fig:Toep:method} 
\end{figure}

We summarize the resulting approach in Algorithm~\ref{alg:repeat}, which has a computational complexity of $\mathcal{O}(n k m^{3k})$, where $m$ is the number of nonzero diagonals in $T$. It proceeds as follows in order to identify and exploit repeated entries along the diagonals of $T^{\ell}$, for $\ell = 0,\ldots,k$: First, the sets $\mathcal{S}_0(T),\ldots,\mathcal{S}_k(T)$ are computed, and then $\mathcal{U}_k(T)$ is constructed (line~1). For each diagonal in $\mathcal{U}_k(T)$, the algorithm determines the starting indices and extracts the initial $\delta$-sets, which are then used to compute the $\delta$-sets of other entries along that diagonal based on Proposition~\ref{prop:facts_delta} (lines~3--5). The diagonal entries are then classified as repeated or non-repeated according to condition~\eqref{cond1} (lines~5--14), and the reduced index set $\Delta_G^{(k)}$ is constructed by collecting the $\delta$-sets corresponding to non-repeated entries. Labeling an entry as non-repeated in line~14 means selecting a representative for all neighboring elements along the diagonal that satisfy condition~\eqref{cond1} with that entry. The matrix function is evaluated only for the submatrix associated with $\Delta_G^{(k)}$ (line~18). Finally, the computed values are assigned to their respective locations in $\widetilde{F}$ and copied to all their associated ``repeated'' positions in order to recover the remaining entries along the diagonals (lines~19--20).

The size of the matrix $X_f$ in \eqref{eq:Xf} is independent of the size of $T$ and actually depends only on the number of the initial diagonals of the matrix and the degree $k$, as the next lemma shows.

\begin{lemma}\label{lem:sub:XF:unchang}
Let $T \in \mathbb{C}^{n \times n}$ be a Toeplitz matrix with $m$ nonzero diagonals, and let $X_f$ be the matrix defined in~\eqref{eq:Xf} for a given $k$. Then, the size of the matrix $X_f$ is of $\mathcal{O}(m^{3k})$.
\end{lemma}
\begin{proof}
It is sufficient to prove the statement for the size of $\Delta_G^{(k)}$. According to~\eqref{Glob:del}, we have
\[
|\Delta_G^{(k)}| \leq |\mathcal{P}| \cdot \max_{(i,j)\in \mathcal{P}} |\Delta_{i,j}^{(k)}|.
\]
From Proposition~\ref{prop:facts_delta} (part 4), we know that the maximum size of any delta-set is less than the size of $\mathcal{U}_k(T)$. In other words, we have
$
|\Delta_G^{(k)}| \leq |\mathcal{P}| \cdot |\mathcal{U}_k(T)|.
$
On the other hand, the condition defining $\mathcal{P}$ in \eqref{P:4:Glob:del} can be rewritten as
\[
\delta_{ij}^{(k)} \cap \{1,\ldots,n\} \neq \left( \left( \delta_{i,j}^{(k)} - 1 \right) \cap \{1,\ldots,n\} \right) + 1.
\]
It is easy to see that this holds when either $1$ or $n+1$ belongs to $\delta_{ij}^{(k)}$. According to the definition of $\delta_{ij}^{(k)}$, when $1$ (and similarly $n+1$) is included in $\delta_{ij}^{(k)}$, the possible cases are
$i = 1 - r_i$ and $j = 1 + r_j,$
where $r_i,r_j\in \mathcal{U}_k(T)$. This shows that the total number of possible pairs $(i,j)$ is at most $|\mathcal{U}_k(T)|^2$. Therefore, we conclude that
$
|\Delta_G^{(k)}| \leq 2|\mathcal{U}_k(T)|^3
$, where the factor $2$ arises from considering both boundaries of the set $\{1,\ldots,n\}$. 
Since $|\mathcal{S}_p(T)|\leq m^p$, we have $|\mathcal{U}_k(T)| = \mathcal{O}(m^k)$, from which the result follows immediately.
\end{proof}

Algorithm~\ref{alg:repeat} relies on an appropriate choice of the polynomial approximation degree $k$. We discuss possible approaches for choosing $k$ in Appendix~\ref{appendix}. The presented approaches  are partly heuristic, and it might happen that they do not yield a satisfactory value of $k$. In that case, it is also possible to run an adaptive version of the algorithm that determines a suitable value of $k$ by increasing the value until a satisfactory accuracy is reached. To do so, we need to be able to estimate the accuracy of the approximation. This can be done by comparing approximations obtained for different degrees $k$ and $\widetilde{k} > k$.

This can be done as follows. A very simple estimate is
\[
\Vert f(T) - \widetilde{F}_{k} \Vert_{\infty} \approx \Vert \widetilde{F}_{\widetilde{k}} - \widetilde{F}_{k} \Vert_{\infty},
\]
where $\widetilde{F}_{k}$ and  $\widetilde{F}_{\widetilde{k}}$ denote the approximations obtained from Algorithm~\ref{alg:repeat} with polynomial degrees $k$ and $\widetilde{k}$, respectively. A satisfactory $k$ is accepted once we have $\Vert \widetilde{F}_{\widetilde{k}} - \widetilde{F}_{k} \Vert_{\infty} < \varepsilon$, for a prescribed threshold value $\varepsilon$. This approach can be made more efficient by only looking at those entries which appear on ``new'' diagonals of $\widetilde{F}_{\widetilde{k}}$ (i.e., those diagonals which were all zero for $\widetilde{F}_k$), as experimentally, these appear to carry most of the weight of the error.

Note that when adaptively increasing the polynomial degree $k$, the computations performed for obtaining $\mathcal{S}_\ell(T), \mathcal{U}_{k}(T)$ and $\Delta_{ij}^{(k)}(T)$ can be reused and one does not need to start the computations from scratch for degree $\widetilde{k}$.


\section{Numerical experiments}\label{sec7}
In this section, we present various numerical experiments on real-world problems and matrices from applications to investigate the quality of our algorithms, approximations, and error bounds. All the experiments were conducted in MATLAB R2021a. To obtain the errors of our approximations (up to machine precision), we compute the exact quantities $f(A)$, $[f(A)]_{ij}$, and $\mathrm{trace}(f(A))$ using the MATLAB commands $\mathtt{expm}$, $\mathtt{sqrtm}$, $\mathtt{inv}$, and $\mathtt{trace}$.

\begin{figure}[t] 
\centering
\begin{subfigure}[b]{0.47\textwidth}
\centering
\includegraphics[width=0.95\textwidth]{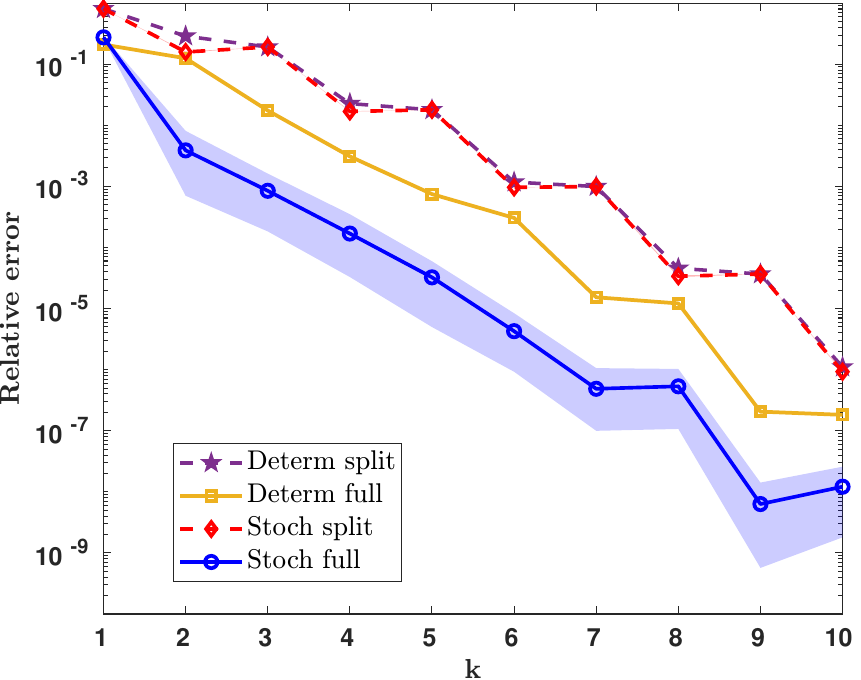}
\caption{\texttt{Gset/G50}, $f(x) = \exp(-x)$.}\label{a-Fig63}
\end{subfigure}
\hfill
\begin{subfigure}[b]{0.47\textwidth}
\centering
\includegraphics[width=0.95\textwidth]{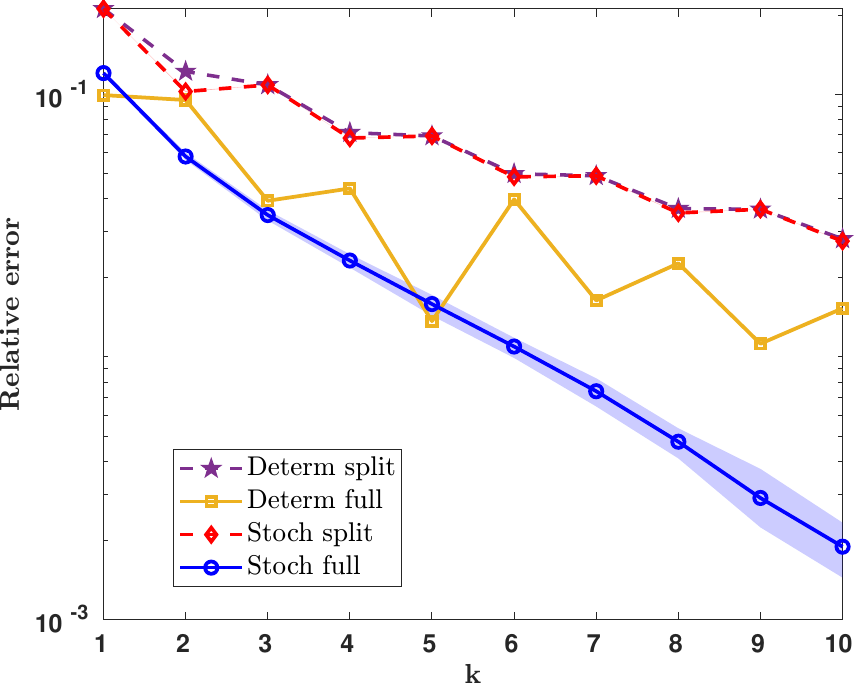}
\caption{\texttt{Gset/G50}, $f(x) = x^{-1/2}$.}
\end{subfigure}

\vspace{.75cm}
\begin{subfigure}[b]{0.47\textwidth}
\centering
\includegraphics[width=0.95\textwidth]{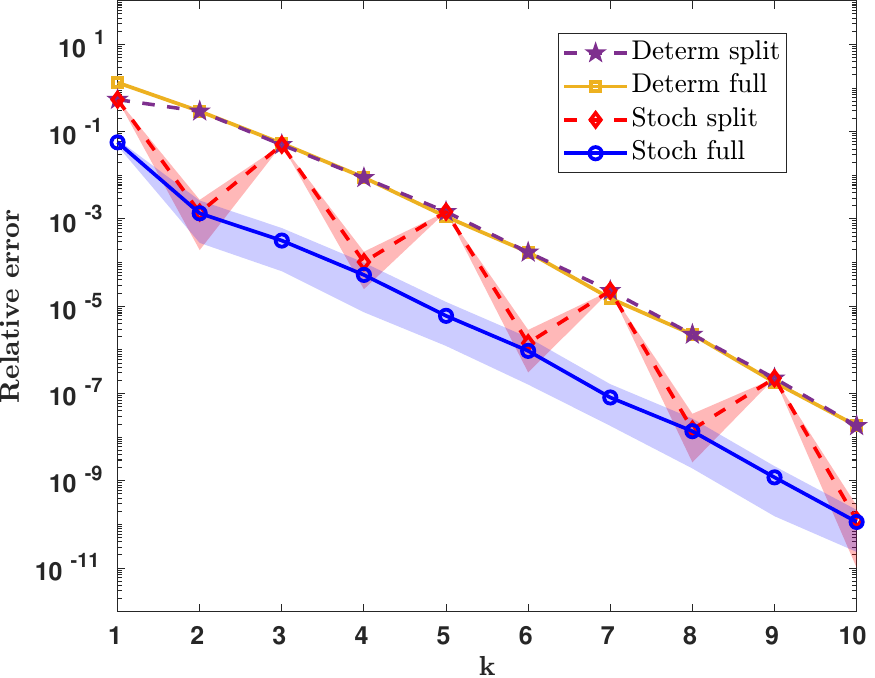}
\caption{\texttt{Norris/fv2} , $f(x) = \exp(-x)$.}\label{c-Fig63}
\end{subfigure}
\hfill
\begin{subfigure}[b]{0.47\textwidth}
\centering
\includegraphics[width=0.95\textwidth]{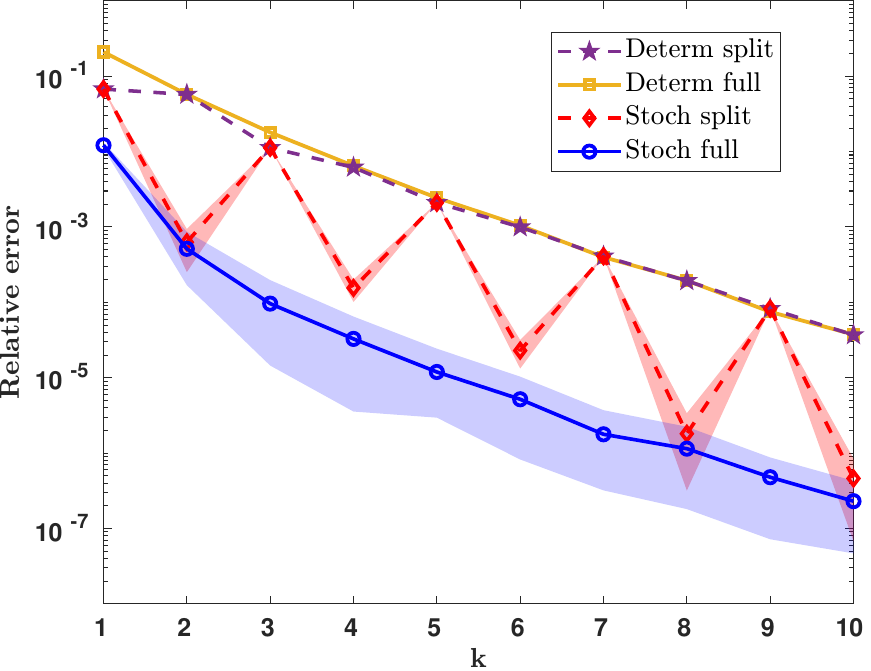}
\caption{\texttt{Norris/fv2} , $f(x) = x^{-1/2}$.}
\end{subfigure}
\caption{Comparison of relative errors of stochastic and deterministic methods for approximating the trace of $f(A)$, for $A$ with $k = 1,\dots,10$, and with $100$ repetitions for each $k$ for the stochastic methods.}\label{ktrace:err:Fig:3}
\end{figure}

\subsection{Scaling of trace approximation methods with $k$}
In this experiment, we illustrate how the trace approximation methods presented in Section~\ref{sec5} behave when increasing the degree $k$ of the polynomial approximation. We use the matrices \texttt{Gset/G50} and \texttt{Norris/fv2} from the \emph{SuiteSparse Matrix Collection}~\cite{davis2011university}, which have sizes $n=3000$ and $n = 9801$, respectively. We consider the functions $f(x) = \exp(x)$ and $f(x) = x^{-1/2}$. As the matrix \texttt{Gset/G50} is indefinite, we replace $A$ by $A+4I$ when computing the inverse square root in order to move all eigenvalues to the positive real line, so that $A^{-1/2}$ is well-defined.

In Figure~\ref{ktrace:err:Fig:3}, we compare the relative errors of the deterministic \texttt{full} and \texttt{split} methods, along with their stochastic counterparts with one vector per set of the partition (i.e., $N_{\ell} \equiv 1$), for $k = 1,\dots,10$. For the stochastic methods, we report the average relative errors for 100 repetitions for each value of $k$.

For the matrix \texttt{Gset/G50}, we observe that the stochastic and deterministic \texttt{split} behave very similarly. For the exponential function, the \texttt{full} deterministic method has an error that is roughly one order of magnitude smaller and the stochastic \texttt{full} method turns to have a smaller error by about one order of magnitude. This is a well-known phenomenon for stochastic probing methods; see~\cite{frommer2025analysis}. The difference between the \texttt{full} and \texttt{split} method as well as the fact that the stochastic \texttt{split} method does not improve upon its deterministic counterpart indicate that the error introduced by the splitting dominates in this case, so that the stochastic improvement in the actual trace estimation does not affect the overall accuracy. For the inverse square root, the improvement of the stochastic over the deterministic \texttt{full} method is more pronounced for most values of $k$.

For \texttt{Norris/fv2}, the results look slightly different. Here, both deterministic methods are on par, indicating that the splitting introduces an error that is negligible. The stochastic \texttt{full} method again consistently has the smallest error. The stochastic \texttt{split} method greatly improves upon its polynomial counterpart only for every other degree, which is likely related to special structural properties of the matrix at hand. Still, it is always at least as accurate as the deterministic method, indicating that it seems always advisable to use the stochastic method. In contrast to the previous example, the results are very similar for both considered functions.

Another interesting observation is that the \texttt{split} method has a much smaller variance than the \texttt{full} method across all values of $k$ and for both functions $f$. This can be explained by the fact that in Theorem \ref{thm:new_bound_variance}, we actually have the series of inequalities
\[
V_{\ell} \leq 2 |P_{\ell}| \min_{p_k\in \Pi_k} \max_{z\in\mathcal{W}(B_{P_{\ell}})}\vert f-p_k\vert^2 \leq 2 |P_{\ell}| \min_{p_k\in \Pi_k} \max_{z\in\mathcal{W}(A)}\vert f-p_k\vert^2,
\]
suggesting that the variance for the \texttt{split} method \emph{could}, in principle, be much smaller than that of the \texttt{full} method.

\subsection{Comparison of algorithms for functions of Toeplitz matrices}\label{Example:Toeplitz_Fun} 
In this experiment, we compare Algorithm~\ref{alg:repeat} for computing functions of Toeplitz matrices with other algorithms from the literature, the $\mathtt{sexpmt}$ method of~\cite{Kressner-Luce} and the methods provided in the CQT toolbox~\cite{bini2019quasi}. We also report results for the built-in MATLAB function $\mathtt{expm}$ as the baseline.
The example we consider is inspired by quantum walks on circulant graphs; see, e.g.,~\cite{todtli2016continuous,yu2024classification}. Specifically, we let $G$ be a graph with $n$ nodes, where each node is connected to its left and right neighbors (and nodes are arranged in a circle, that is, node $n$ is connected to node $1$). Additionally, there are connections at a specified distance $g > 1$. For example, when $g = 5$, node 1 is connected to node 6, node 2 is connected to node 7, and so on, until node $n-1$ is connected to node 4, and node $n$ is connected to node 5. We let $A = L_G$, the Laplacian of the graph $G$, so that for a time step $t > 0$, the matrix function $\exp(tA)$ is the heat diffusion propagator on the graph.

\begin{table}[htb]
\caption{Run times and relative errors (using the infinity norm) for computing $\exp(tA)$, where $A$ is a graph Laplacian. Note that for $n = 25,\!000$, calling \texttt{expm} required more than the 16 GB main memory that was available.}\label{tab:Top:comp}
\centering
\begin{tabular}{c||ccc|cc|cc|cc|c}
\toprule
&
& 
&
& \multicolumn{2}{c|}{Algorithm 2}
& \multicolumn{2}{c|}{\texttt{CQT}}
& \multicolumn{2}{c|}{\texttt{sexpmt}}
& \texttt{expm} \\
$g$ & $n$ & $k$ & $|\Delta|$
 & Time & Err.
 & Time & Err.
 & Time & Err.
 & Time \\
\midrule
{} & 1000 & 6 & 25 & 0.05 & 9.1e-12 & 0.06 & 2.1e-12 & 0.69 & 2.3e-11 & 0.08 \\
{} & 5000 & 6  & 25 & 0.21 & 3.2e-13 & 0.54 & 2.2e-12 & 10.34 & 9.1e-11 & 11.59 \\
$2$ & 10000 & 6 & 25 & 0.51 & 1.3e-12 & 1.60 & 3.1e-12 & 37.04 & 2.1e-10 & 77.97 \\
{} & 15000 & 6 & 25 & 0.92 & 1.4e-12 & 3.74 & 3.7e-12 & 81.38 & 3.7e-10 & 255.90 \\
{} & 20000 & 6 & 25 & 2.08 & 2.5e-12 & 10.19 & 3.1e-12 & 154.76 & 4.1e-10 & 660.16 \\
{} & 25000 & 6 & 25 & 5.08 &  & 26.34 &  & 235.45 &  &  \\
\midrule
{} & 1000  & 6 & 53 & 0.10 & 2.1e-13 & 0.10 & 2.1e-12 & 0.70 & 3.2e-11 & 0.07 \\
{} & 5000  & 6 & 53 & 0.42 & 4.2e-13 & 0.82 & 2.2e-12 & 10.90 & 2.2e-10 & 11.61 \\
$5$ & 10000 & 6 & 53 & 0.95 & 1.4e-12 & 2.12 & 3.1e-12 & 39.29 & 4.3e-10 & 73.64 \\
{} & 15000 & 6 & 53 & 1.58 & 2.1e-12 & 4.70 & 3.7e-12 & 85.77 & 6.5e-10 & 260.90 \\
{} & 20000 & 6 & 53 & 2.59 & 2.3e-12 & 9.73 & 4.1e-12 & 142.67 & 8.1e-10 & 657.30 \\
{} & 25000 & 6 & 53 & 6.08 &  & 29.04 &  & 212.36 &  &  \\
\midrule
{} & 1000  & 6 & 85 & 0.31 & 4.5e-13 & 1.14 & 2.3e-12 & 1.27 & 3.1e-11 & 0.09 \\
{} & 5000  & 6 & 85 & 0.62 & 6.7e-13 & 2.93 & 2.8e-12 & 9.91 & 2.2e-10 & 11.33 \\
$20$ & 10000 & 6 & 85 & 1.34 & 1.8e-12 & 6.59 & 3.7e-12 & 36.01 & 3.6e-10 & 68.08 \\
{} & 15000 & 6 & 85 & 2.31 & 2.1e-12 & 11.68 & 3.1e-12 & 79.63 & 6.5e-10 & 246.93 \\
{} & 20000 & 6 & 85 & 6.01 & 2.1e-12 & 24.40 & 4.4e-12 & 141.81 & 8.1e-10 & 635.35 \\
{} & 25000 & 6 & 85 & 7.27 &  & 41.58 &  & 227.04 &  &  \\
\bottomrule
\end{tabular}
\end{table}

In our experiment, we take the time step $t = 0.01$, coupling distances $g\in \{2, 5, 20\}$ and varying graph sizes $n$ ranging from $1000$ to $25,000$. To allow for a fair comparison, we perform a preprocessing step before using the CQT toolbox: its methods work most efficiently when representing $A = T+E$, where $T$ is a banded Toeplitz matrix and $E$ is of low-rank. Clearly, one can find such a representation for $A$, where the rank of the correction term $E$ is $2g$. Rewriting the problem in this way greatly speeds up the CQT methods, and we do not include the time for the preprocessing step in the timings we report. 

Our results are summarized in Table~\ref{tab:Top:comp}. For all considered methods, we report the run time and the final approximation error (compared to the solution computed by the baseline method \texttt{expm}). We also provide the polynomial approximation degree $k$ (which is automatically determined using the approach outlined in Appendix~\ref{appendix}) and the corresponding submatrix size, $|\Delta|$, used in Algorithm~\ref{alg:repeat}. As expected from Lemma~\ref{lem:sub:XF:unchang}, the size of the submatrix does not increase when the matrix size $n$ is increased.

We observe that, in all cases, our algorithm outperforms the competing methods. In particular, it appears to scale more favorably with the matrix size $n$.

%
\subsection{Estrada index}\label{Example:Trace_Approx} 
For a given graph $G$, the Estrada index of $G$ is defined as $\mathrm{trace}(\exp(A_G))$, where $A_G$ is the adjacency matrix of the graph $G$. This index measures the strength of connectivity within $G$~\cite{Estrada2012}.

In our experiment, we aim to approximate the Estrada index of graphs from the SuiteSparse collection using the stochastic delta estimator, both in the \texttt{full} and \texttt{split} versions, and compare them with other established trace estimation methods, namely XTrace~\cite{epperly2024xtrace} and A-Hutch++ \cite{meyer2021hutch++,persson2022improved}.

For our method, given a tolerance $\varepsilon>0$, we computed candidate values for the polynomial approximation degree $k$ by the methodology outlined in \ref{appendix}.\footnote{With $\epsilon = \mathrm{trace}(\exp(A))\varepsilon/2n$, $\tau = 1.1$ and $\eta = 0.05, 0.06, 0.07, \dots, 1$. For large matrices, we used $\mathrm{trace}(\exp(A)) \geq n + m$, where $m$ is the number of edges.} To find the value of $k$ to be used, we use the stopping criterion,
\[
 \dfrac{\left\vert \mathcal{RT}_{\Delta}^{(k_1)}(\exp(A)) - \mathcal{RT}_{\Delta}^{(k_2)}(\exp(A))\right\vert}{\left\vert\mathcal{RT}_{\Delta}^{(k_2)}(\exp(A))\right\vert}< \varepsilon,
\] 
based on two consecutive candidate values $k_1<k_2$. This stopping criterion turned out to be rather conservative and may actually provide approximations that are more accurate than desired. We set the number of vectors for each partition to $N_\ell = 1$.

In all methods, we use the Lanczos method to approximate matrix multiplication with $\exp(A)$. For a fair run time comparison, we use the same number of Lanczos iterations in all cases. For our delta trace estimators, we follow the recommendation for probing methods of~\cite{Frommer-Schimmel-Schweitzer} to use $2k$ iterations, and then use the same number in A-Hutch++ and XTrace. We observed that using this number of iterations, the Lanczos approximation error becomes negligible compared to the trace estimation error in these methods. 
Additionally, in the A-Hutch++ method, we set the failure probability to $\delta = 0.05$.\footnote{We used \texttt{xtrace\_tol} from \url{https://github.com/eepperly/XTrace/tree/main} for \texttt{XTrace}, and \texttt{adap\_hpp} from \url{https://github.com/davpersson/A-Hutch-} for \texttt{A-Hutch++}.}

We set the tolerance to $\varepsilon = \mathtt{1e-4}$ for all methods and report the resulting computing times for all the considered methods in Table \ref{estrada:new}, along with the relative error of the methods.\footnote{For large matrices, instead of comparing with \texttt{expm}, we compare our method with the output of \texttt{xtrace\_tol} to show the accuracy of our method.} We observe that our methods are very competitive as compared to the established algorithms, and, for all but two problems, one of our two methods is always the fastest. We also note that for many problems, our methods yield approximations that are many orders of magnitude smaller than requested, but are at the same time still much faster than the competing methods with larger errors.

\begin{table}[htb]
\caption{Runtime and relative error in computing the Estrada index of some graphs  from the \emph{SuiteSparse} collection using different stochastic methods. For each problem, the run time of the fastest method is marked in bold.}\label{estrada:new}
\centering
\scalebox{.9}{%
\begin{tabular}{c|c c|cc|cc|cc|cc|c}
\toprule
&
&
& \multicolumn{2}{c|}{\texttt{full}} 
& \multicolumn{2}{c|}{\texttt{split}} 
& \multicolumn{2}{c|}{\texttt{XTrace}} 
& \multicolumn{2}{c|}{\texttt{A-Hutch++}} 
& \texttt{expm}\\
\texttt{ID} & $n$  & $k$
 &  Time & Err.
 & Time & Err.
 & Time & Err.
 & Time & Err.
 &   Time\\ 
\midrule
\texttt{139} & 992 &  11 & 0.48 & 6.4e-05 & 0.70 & 1.6e-04 & 0.08 & 1.0e-05 & 0.09 & 2.9e-05 & \textbf{0.07} \\ 
\texttt{511} & 3000 & 5 & 0.14 & 3.5e-05 & \textbf{0.10} & 6.2e-06 & 4.24 & 1.7e-05 & 28.78 & 2.0e-05 & 2.10 \\ 
\texttt{512} & 3000 & 6 & 0.20 & 4.7e-06 & \textbf{0.13} & 2.4e-05 & 4.45 & 8.7e-07 & 29.28 & 3.2e-07 & 2.20 \\ 
\texttt{514} & 3000 & 6 & \textbf{0.11} & 2.1e-05 & 0.14 & 1.8e-06 & 4.17 & 3.7e-05 & 28.93 & 1.3e-05 & 2.06 \\ 
\texttt{250} & 3,345 & 5 & 0.73 & 5.2e-07 & \textbf{0.72} & 9.8e-09 & 0.78 & 1.2e-03 & 3.09 & 4.2e-05 & 3.10 \\ 
2411 & 4,253 & 4 & 4.63 & 1.3e-12 & 4.69 & 1.3e-12 & \textbf{2.26} & 8.1e-05 & 23.03 & 2.6e-05 & 7.34 \\ 
\texttt{2412} & 8,034 & 4 & 10.82 & 6.6e-15 & \textbf{10.81} & 6.6e-15 & 26.53 & 1.4e-04 & 445.97 & 2.0e-06 & 35.83 \\ 
\texttt{2527} & 16,386 & 5 & \textbf{6.69} & 6.1e-03 & 6.77 & 6.1e-03 & 66.62 & 6.1e-03 & 605.35 & 6.0e-03 & 361.37 \\ 
\texttt{1255} & 17,222 & 8 &16.58 & 1.4e-04 & \textbf{16.53} & 1.9e-04 & 17.12 & 6.8e-06 & 43.03 & 2.8e-05 & 383.98 \\ 
\texttt{2438} & 49,152 & 4 & \textbf{2.93} & 9.8e-05  & 8.97 & 6.3e-05  & 489.3 & $\ast$  &  $>600$ & ---  & $>600$ \\ 
\texttt{2437} & 49,152 & 4 & \textbf{7.37} & 9.3e-05  & 7.41 & 9.5e-05  & 500.7 & $\ast$  &  $>600$ & ---  & $>600$ \\ 
\texttt{2435} & 49,152 & 5 & 69.62 & 5.5e-05  & \textbf{69.44} & 7.9e-05  & 493.4 & $\ast$  &  $>600$ & ---  & $>600$ \\ 
\bottomrule
\end{tabular}}%
\end{table}
\section{Concluding remarks}\label{sec8}
We introduced a new methodology for identifying the locations where nonzero entries in matrix polynomials $p(A)$ can occur, which worked particularly well when $A$ is a matrix in which not too many sub- or superdiagonals contain nonzero entries. 

Our methodology can be used very efficiently in the context of probing methods for estimating the trace of matrix functions. These methods require a partitioning of the rows and columns of $A$, which is typically determined by computing a distance-coloring of the graph of $A$, a step that can be extremely costly for large matrices. In this context, our techniques can be used to obtain an alternative preprocessing procedure that is often much more efficient. Additionally, it can help to further speed up probing methods by reducing the size of the subproblems  needed to be solved.

As a second application of our methodology, we discussed how we could use it to compute functions of Toeplitz matrices $T$ very efficiently: it allowed us to find a small submatrix of $T$ that essentially determined the whole matrix polynomial $p(T)$, which could, in turn, be used to approximate $f(T)$. 

We illustrated in numerical experiments that, in situations where it can be efficiently applied, our method is very competitive with state of the art methods for both considered tasks.

An open topic for future research is the development of accurate and reliable error estimates and, based on that, an efficient methodology for determining a suitable polynomial approximation degree to be used in our algorithms. 
\section{Acknowledgements}
Our sincere thanks are extended to Daniel Kressner for providing the MATLAB codes for \texttt{expmt} and to Dario Bini for sharing information about the \texttt{CQT toolbox}. We also thank two anonymous reviewers for their careful reading of the manuscript and their valuable comments which helped us a lot to improve the quality of exposition. The second author thanks Sharif University of Technology for its support.

\bibliographystyle{siam}
\bibliography{lit}

@article{Beckermann,
  author  = {Beckermann, B.},
  title   = {Image numérique, {GMRES} et polynomes de {F}aber},
  journal = {C.~R. Acad. Sci. Paris Ser. I},
  volume  = {340},
  year    = {2005},
  pages   = {855--860},
}

@article{Bekas-Kokiopoulou-Saad,
  author  = {Bekas, C. and Kokiopoulou, E. and Saad, Y.},
  title   = {An estimator for the diagonal of a matrix},
  journal = {Appl. Numer. Math.},
  volume  = {57},
  year    = {2007},
  pages   = {1214--1229},
}

@article{Ben_Boit,
  author  = {Benzi, M. and Boito, P.},
  title   = {Decay properties for functions of matrices over {$C^{\ast}$}-algebras},
  journal = {Linear Algebra Appl.},
  volume  = {456},
  year    = {2014},
  pages   = {174--198},
}

@article{Benzi-Boito-Razouk,
  author  = {Benzi, M. and Boito, P. and Razouk, N.},
  title   = {Decay properties of spectral projectors with applications to electronic structure},
  journal = {SIAM Rev.},
  volume  = {55},
  year    = {2013},
  pages   = {3--64},
}

@article{Benzi-Golub,
  author  = {Benzi, M. and Golub, G.~H.},
  title   = {Bounds for the entries of matrix functions with applications to preconditioning},
  journal = {BIT},
  volume  = {39},
  year    = {1999},
  pages   = {417--438},
}

@article{Ben_Raz,
  author  = {Benzi, M. and Razouk, N.},
  title   = {Decay bounds and $\mathcal{O}(n)$ algorithms for approximating functions of sparse matrices},
  journal = {Electron. Trans. Numer. Anal.},
  volume  = {28},
  year    = {2007},
  pages   = {16--39},
}

@article{crouzeix-palencia,
  author  = {Crouzeix, M. and Palencia, C.},
  title   = {The numerical range is a $(1 + \sqrt{2})$-spectral set},
  journal = {SIAM J. Matrix Anal. Appl.},
  volume  = {38},
  year    = {2017},
  pages   = {649--655},
}

@article{Demko,
  author  = {Demko, S. and Moss, W.~F. and Smith, P.~W.},
  title   = {Decay rates for inverses of banded matrices},
  journal = {Math. Comp.},
  volume  = {43},
  year    = {1984},
  pages   = {491--499},
}

@article{Dur,
  author  = {Druskin, V. and G{\"u}ttel, S. and Knizhnerman, L.},
  title   = {Near-optimal perfectly matched layers for indefinite {H}elmholtz problems},
  journal = {SIAM Rev.},
  volume  = {58},
  year    = {2016},
  pages   = {90--116},
}

@article{Estrada-Higham1,
  author  = {Estrada, E. and Higham, D.~J.},
  title   = {Network properties revealed through matrix functions},
  journal = {SIAM Rev.},
  volume  = {52},
  year    = {2010},
  pages   = {696--714},
}

@article{Frommer-Schimmel-Schweitzer,
  author  = {Frommer, A. and Schimmel, C. and Schweitzer, M.},
  title   = {Analysis of probing techniques for sparse approximation and trace estimation of decaying matrix functions},
  journal = {SIAM J. Matrix Anal. Appl.},
  volume  = {42},
  number  = {3},
  year    = {2021},
  pages   = {1290--1318},
}

@article{Goedecker,
  author  = {Goedecker, S.},
  title   = {Linear scaling electronic structure methods},
  journal = {Rev. Modern Phys.},
  volume  = {71},
  year    = {1999},
  pages   = {1085--1123},
}

@book{NHigham,
  author    = {Higham, N. J.},
  title     = {Function of Matrices: {T}heory and Computation},
  publisher = {SIAM},
  address   = {Philadelphia},
  year      = {2008},
}

@article{Kressner-Luce,
  author  = {Kressner, D. and Luce, R.},
  title   = {Fast computation of the matrix exponential for a {T}oeplitz matrix},
  journal = {SIAM J. Matrix Anal. Appl.},
  volume  = {39},
  year    = {2018},
  pages   = {23--47},
}

@article{Lee,
  author  = {Lee, S.~T. and Pang, H.-K. and Sun, H.-W.},
  title   = {Shift-invert {A}rnoldi approximation to the {T}oeplitz matrix exponential},
  journal = {SIAM J. Sci. Comput.},
  volume  = {32},
  year    = {2010},
  pages   = {774--792},
}

@book{Suetin,
  author    = {Suetin, P.~K.},
  title     = {Series of {F}aber polynomials},
  publisher = {Gordon and Breach Science Publishers},
  address   = {Amsterdam},
  year      = {1998},
  note      = {Translated from the 1984 Russian original by E.\,V.\ Pankratiev},
}

@article{FrommerSchimmelSchweitzer2018a,
 author = {Frommer, Andreas and Schimmel, Claudia and Schweitzer, Marcel},
 title = {Bounds for the decay of the entries in inverses and {C}auchy--{S}tieltjes functions of certain sparse, normal matrices},
 journal = {Numer. Linear Algebra Appl.},
 volume = {25},
 number = {4},
 pages = {e2131, 17},
 year = {2018},
 doi = {10.1002/nla.2131},
}

@article{tang2012probing,
  title={A probing method for computing the diagonal of a matrix inverse},
  author={Tang, Jok M and Saad, Yousef},
  journal={Numer. Linear Algebra Appl.},
  volume={19},
  number={3},
  pages={485--501},
  year={2012},
  publisher={Wiley Online Library}
}

@article{frommer2025analysis,
  title={Analysis of stochastic probing methods for estimating the trace of functions of sparse symmetric matrices},
  author={Frommer, Andreas and Rinelli, Michele and Schweitzer, Marcel},
  journal={Math. Comp.},
  volume={94},
  number={352},
  pages={801--823},
  year={2025}
}

@article{stathopoulos2013hierarchical,
  title={Hierarchical probing for estimating the trace of the matrix inverse on toroidal lattices},
  author={Stathopoulos, Andreas and Laeuchli, Jesse and Orginos, Kostas},
  journal={SIAM J. Sci. Comput.},
  volume={35},
  number={5},
  pages={S299--S322},
  year={2013},
  publisher={SIAM}
}

@article{benzi2023computation,
  title={Computation of the von {N}eumann entropy of large matrices via trace estimators and rational {K}rylov methods},
  author={Benzi, Michele and Rinelli, Michele and Simunec, Igor},
  journal={Numer. Math.},
  volume={155},
  number={3},
  pages={377--414},
  year={2023},
  publisher={Springer}
}

@Article{bere09,
  author = "Beckermann, Bernhard and Reichel, Lothar",
  title = "Error Estimates and Evaluation of Matrix Functions via the {F}aber
                  Transform",
  journal = "SIAM J. Numer. Anal.",
  year = 2009,
  volume = 47,
  number = 5,
  pages = "3849–3883",
  month = jan,
  issn = "1095-7170",
  doi = "10.1137/080741744"
  }

@Article{drkn89,
  author = "Druskin, V.L. and Knizhnerman, L.A.",
  title = "Two polynomial methods of calculating functions of symmetric matrices",
  journal = "USSR Computat. Math. Math. Phys.",
  year = 1989,
  volume = 29,
  number = 6,
  pages = "112–121",
  month = jan,
  issn = "0041-5553",
  doi = "10.1016/s0041-5553(89)80020-5"
  }

@article{Saad1992,
	Author = {Yousef Saad},
	journal = {SIAM J. Numer. Anal.},
	Month = {February},
	Number = {1},
	Pages = {209--228},
	Title = {Analysis of Some {K}rylov Subspace Approximations to the matrix Exponential Operator},
	Volume = {29},
	Year = {1992}}

@phdthesis{Guettel2010,
	Author = {Stefan G\"{u}ttel},
	School = {{Fakult\"{a}t f\"{u}r Mathe\-matik und Informatik der Technischen Universit\"{a}t Bergakademie Freiberg}},
	Title = {Rational {K}rylov Methods for Operator Functions},
	Type = {{PhD thesis}},
	Year = {2010},
}

@phdthesis{schimmel2019bounds,
  title={Bounds for the decay in matrix functions and its exploitation in matrix computations},
  author={Schimmel, Claudia},
  year={2019},
  school={Dissertation, Wuppertal, Bergische Universit{\"a}t, 2019}
}

@article{Barnes1989,
  author  = {E. R. Barnes},
  title   = {Circular Discs Containing Eigenvalues of Normal Matrices},
  journal = {Linear Algebra Appl.},
  volume  = {114/115},
  year    = {1989},
  pages   = {501--521}
}

@article{BhatiaSharma2012,
  author  = {R. Bhatia and R. Sharma},
  title   = {Some inequalities for positive linear maps},
  journal = {Linear Algebra Appl.},
  volume  = {436},
  year    = {2012},
  pages   = {1562--1571}
}

@article{KovariPommerenke1967,
  author  = {T. Kővari and Ch. Pommerenke},
  title   = {On {F}aber polynomials and {F}aber expansions},
  journal = {Math. Z.},
  volume  = {99},
  year    = {1967},
  pages   = {193--206}
}

@article{Merikoski2003,
  author  = {J. K. Merikoski and R. Kumar},
  title   = {Characterisations and lower bounds for the spread of a normal matrix},
  journal = {Linear Algebra Appl.},
  volume  = {364},
  year    = {2003},
  pages   = {13--31}
}

@article{Mirsky1956,
  author  = {L. Mirsky},
  title   = {The spread of a matrix},
  journal = {Mathematika},
  volume  = {3},
  year    = {1956},
  pages   = {127--130}
}

@article{TrefethenWeideman2014,
  author  = {L. N. Trefethen and J. A. C. Weideman},
  title   = {The Exponentially Convergent Trapezoidal Rule},
  journal = {SIAM Rev.},
  volume  = {56},
  year    = {2014},
  pages   = {385--458}
}

@article{Hutch,
  title={A stochastic estimator of the trace of the influence matrix for {L}aplacian smoothing splines},
  author={Hutchinson, Michael F},
  journal={Commun.\ Stat.\ Simul.\ Comput.},
  volume={18},
  number={3},
  pages={1059--1076},
  year={1989},
  publisher={Taylor \& Francis}
}

@book{golub2009matrices,
  title={Matrices, Moments and Quadrature with Applications},
  author={Golub, Gene H and Meurant, G{\'e}rard},
  year={2009},
  publisher={Princeton University Press}
}

@article{yu2024classification,
  title={On the classification and dispersability of circulant graphs with two jump lengths},
  author={Yu, Xiaoxiang and Shao, Zeling and Li, Zhiguo},
  journal={Discrete Appl.\ Math.},
  volume={355},
  pages={268--286},
  year={2024},
  publisher={Elsevier}
}

@article{artalejo2010markovian,
  title={Markovian arrivals in stochastic modelling: a survey and some new results},
  author={Artalejo, Jes{\'u}s R and G{\'o}mez-Corral, Antonio},
  journal={SORT},
  volume={34},
  number={2},
  pages={101--156},
  year={2010}
}

@article{todtli2016continuous,
  title={Continuous-time quantum walks on directed bipartite graphs},
  author={T{\"o}dtli, Beat and Laner, Monika and Semenov, Jouri and Paoli, Beatrice and Blattner, Marcel and Kunegis, J{\'e}r{\^o}me},
  journal={Phys.\ Rev.\ A},
  volume={94},
  number={5},
  pages={052338},
  year={2016},
  publisher={APS}
}

@article{davis2011university,
  title={The {U}niversity of {F}lorida sparse matrix collection},
  author={Davis, Timothy A and Hu, Yifan},
  journal={ACM Trans.\ Math.\ Softw.},
  volume={38},
  number={1},
  pages={1--25},
  year={2011},
  publisher={ACM New York, NY, USA}
}

@article{bini2019quasi,
  title={{Q}uasi-{T}oeplitz matrix arithmetic: a {MATLAB} toolbox},
  author={Bini, Dario A and Massei, Stefano and Robol, Leonardo},
  journal={Numer.\ Algorithms},
  volume={81},
  number={2},
  pages={741--769},
  year={2019},
  publisher={Springer}
}

@book{Estrada2012,
  title={The Structure of Complex Networks: Theory and Applications},
  author={Estrada, Ernesto},
  year={2012},
  publisher={Oxford University Press}
}

@article{epperly2024xtrace,
  title={Xtrace: {M}aking the most of every sample in stochastic trace estimation},
  author={Epperly, Ethan N and Tropp, Joel A and Webber, Robert J},
  journal={SIAM J.\ Matrix Anal.\ Appl.},
  volume={45},
  number={1},
  pages={1--23},
  year={2024},
  publisher={SIAM}
}

@inproceedings{meyer2021hutch++,
  title={Hutch++: {}Optimal stochastic trace estimation},
  author={Meyer, Raphael A and Musco, Cameron and Musco, Christopher and Woodruff, David P},
  booktitle={Symposium on Simplicity in Algorithms (SOSA)},
  pages={142--155},
  year={2021},
  organization={SIAM}
}

@article{persson2022improved,
  title={Improved variants of the {H}utch++ algorithm for trace estimation},
  author={Persson, David and Cortinovis, Alice and Kressner, Daniel},
  journal={SIAM J.\ Matrix Anal.\ Appl.},
  volume={43},
  number={3},
  pages={1162--1185},
  year={2022},
  publisher={SIAM}
}

@article{stewart2002,
  title={A {K}rylov--{S}chur algorithm for large eigenproblems},
  author={Stewart, Gilbert W},
  journal={SIAM J.\ Matrix Anal.\ Appl.},
  volume={23},
  number={3},
  pages={601--614},
  year={2002},
  publisher={SIAM}
}

@article{Gray2006,
  title={Toeplitz and {C}irculant {M}atrices: {A} review},
  author={Gray, Robert M},
  year={2006},
  publisher={now publishers inc}
}
\appendix
\section{Heuristic method for finding a suitable polynomial degree}\label{appendix}
Our methods (as well as many other matrix function algorithms based on polynomial approximation) crucially rely on determining a suitable polynomial degree $k$ such that a specified accuracy is reached. Here, we brief\/ly sketch an effective (heuristic) approach that yields \emph{candidate values} of $k$ that can then be tested, along with a suitable error estimate (cf.~, e.g., the discussion at the end of Section~\ref{sec6}).

We begin by recalling some basic facts about Faber polynomials; cf.~\cite{Beckermann,bere09,drkn89}. A set $\mathcal{K} \subseteq \mathbb{C}$ is called a \emph{continuum} if $\mathcal{K}$ is compact, connected, and contains more than a single point. If $\mathcal{K}$ has a connected complement, then the Riemann mapping theorem ensures the existence of a function $\phi$ that maps the exterior of $\mathcal{K}$ conformally onto the set $\lbrace z \in \mathbb{C} : |z| > 1\rbrace$ so that $\phi(\infty) = \infty,$ and $ \lim_{z\rightarrow \infty} \phi(z)/z = d>0$. The mapping $\phi$ has a Laurent expansion, $\phi(z) = dz + \sum_{j=0}^{\infty} \dfrac{a_j}{z^j}$, so that for every $k \geq 0$, we have
\[
\big(\phi(z)\big)^k = d^kz^k + a_{k-1}^{(k)}z^{k-1}+\cdots+a_0^{(k)}+\sum_{j=1}^{\infty}\dfrac{a_{-j}^{(k)}}{z^j}.
\]
The polynomial parts, $F_k(z) = d^k z^k + a_{k-1}^{(k)}z^{k-1}+\cdots+a_0^{(k)}$, are known as the \textit{Faber polynomials} with respect to $\mathcal{K}$. The following theorem of~\cite{Suetin} shows that any analytic function can be expanded in a Faber series.

\begin{theorem}\label{faber1}
Every function $f(x)$ analytic on a continuum $\mathcal{K}$ can be expanded in a Faber series converging uniformly on $\mathcal{K}$, that is, for any $x\in\mathcal{K}$,
\begin{equation}\label{coeff}
f(x) =\sum_{k=0}^{\infty} a_k F_k(x),~~~~~~~~\text{with}~~ a_k := \dfrac{1}{2\pi {\bf i}}\int_{|z| = \tau } \dfrac{f(\phi^{-1}(z))}{z^{k+1}} ~\mathrm{d}z,
\end{equation}
 where $\tau>1$ is selected such that $f$ is analytic on the complement of {\small $\lbrace \phi^{-1}(z): |z| > \tau \rbrace$}.
\end{theorem}

By~\cite[Theorem 2]{KovariPommerenke1967}, we have
\[
\min_{p_k \in \Pi_k} \max_{z \in \mathcal{K}} |f(z) - p_k(z)| \leq 2 \sum_{m=k+1}^{\infty} |a_m|,
\]
showing that the Faber coefficients provide information on the best possible polynomial approximation error. If we choose the continuum $\mathcal{K}$ such that $\mathcal{W}(A) \subseteq \mathcal{K}$, it is then sufficient to choose $k$ such that the tail of the Faber series satisfies $\sum_{m=k+1}^{\infty} |a_m| < \frac{\epsilon}{2}$ to obtain a polynomial approximation accuracy $\epsilon$ on $\mathcal{W}(A)$. This approach requires computing or estimating the Faber coefficients. In general, one can approximate the coefficients by applying the trapezoidal rule to the integral representation in~\eqref{coeff}, which is known to converge exponentially~\cite{TrefethenWeideman2014}.

A challenge in using Faber coefficients for determining $k$ is finding a continuum $\mathcal{K}$ that encloses the field of values of $A$, which is in general unknown and very costly to compute, and at the same time exhibits a manageable representation of the function $\phi(z)$. A simple, general purpose approach for this task is to use an enclosing circle for $\mathcal{W}(A)$, which can often be computed at a reasonable cost; see, e.g.,~\cite{Mirsky1956, Barnes1989, BhatiaSharma2012, Merikoski2003}. Therefore, in the following we restrict our discussion to the case that $\mathcal{K}$ is a circle with center $c$ and radius $R$, so that $\phi(z) = \frac{z - c}{R}$ and $\phi^{-1}(z) = Rz + c$ and the $k$th Faber coefficient is simply the $k$th Taylor coefficient of $f$ multiplied by $R^k$.

One can obtain such a circle by first finding the smallest possible intervals containing the spectrum of the Hermitian and the skew-Hermitian part of $A$, respectively, denoted by $I_H$ and $I_S$. Then, the circle $\mathcal{K}$  can be constructed using these intervals by setting the radius $R = \sqrt{\ell_H^2+\ell_S^2}/2$, where $\ell_H$ and $\ell_S$ are the lengths of the respective intervals, and the center of the circle is $c = c_{H} + \textbf{i} c_S$, where $c_H$ and $c_S$ are the centers of the corresponding intervals. One efficient way to find such intervals for real sparse matrices is to use MATLAB command \texttt{eigs} for computing only the largest and the smallest eigenvalues of each part \cite{stewart2002}. 
 
Another simple way to find such intervals is by applying the Gershgorin disk theorem. This is particularly efficient for Toeplitz matrices (which we focused on in Section~\ref{sec6}): in this case, all Gershgorin disks have the same center, so that their union is simply the circle with the largest radius. Denoting by $Q$ either the Hermitian or the skew-Hermitian part of the Toeplitz matrix, the length of the intervals can be explicitly determined and is given by $\frac{\ell_Q}{2} = \|Q\|_{\infty} - |[Q]_{11}|$, while the center of the interval corresponding to each part is $c_Q = [Q]_{11}$. We add that the infinity norm of a Toeplitz matrix, as being required for determining the interval lengths, can be computed with maximum complexity $\mathcal{O}(n)$~\cite[Lemma 4.1]{Kressner-Luce}.  

One shortcoming of all methods working with a superset $\mathcal{K} \supset \mathcal{W}(A)$ is that even when $f$ is analytic on $\mathcal{W}(A)$, the enclosing continuum $\mathcal{K}$ might contain a singularity or branch cut of $f$, leading to a deterioration in polynomial approximation quality. This is no problem for entire functions like the exponential, but might be relevant for other functions. It can, therefore, be reasonable to try to shrink the enclosing circle. To this end, we define a set of reduction ratios $ \{ \eta_1,~ \eta_2, \ldots, \eta_p \}$ with $\eta_p = 1$ and $0 < \eta_j<\eta_{j+1}$, for $j = 1,\ldots,p-1$, and consider the sequence of circles with center $c$ and radii $\eta_jR, j = 1,\dots,p$. For each of these circles, the corresponding Faber coefficients~\eqref{coeff} can be computed in order to obtain a corresponding candidate value for $k$. These values can then be tested one after the other (from small to large) until the desired accuracy is reached. This way, the set of possible polynomial degrees to consider is reduced to a small, finite number of values.\footnote{We mention that it is not necessary to include a ``safety mechanism'' to ensure that the ``shrunk circles'' still enclose $\mathcal{W}(A)$. If they do not, then this will be reflected in poor accuracy of the computed approximation---which will be detected by the error estimate---and the corresponding value for $k$ will be discarded.} For appropriate choices of $\eta_j$, in our experiments the first or second candidate value was typically already sufficient to reach the prescribed accuracy.\footnote{In numerical experiments, we considered the case  $ \{ \eta_1,~ \eta_2, \ldots, \eta_p \} = \{0.05, 0.06, 0.07, \dots, 1\}$.}

\end{document}